\documentclass[3p]{my-elsarticle}

\usepackage{floatrow}
\usepackage{hyperref}
\usepackage{booktabs}
\usepackage{array}
\usepackage{graphicx}
\usepackage[label font=bf]{subfig}
\usepackage{float}
\usepackage{placeins}
\usepackage{rotating}
\usepackage{tablefootnote}
\usepackage{threeparttable}
\usepackage{multirow}
\usepackage[tableposition = top]{caption}
\usepackage{amsmath}
\usepackage{amsfonts}
\usepackage{algorithm}
\usepackage{algpseudocode}
\usepackage{tabularx}
\usepackage{pdflscape} 
\usepackage{longtable}
\usepackage{soul}
\usepackage{acro}[=v3]
\usepackage{supertabular}
\acsetup{
    first-style = long-short,
    list/display = used,
    list/template = tabular,
}

\DeclareAcronym{GSA}{ 
    short = {GSA}, 
    long  = {Gravitational Search Algorithm},
    tag = {abbrev}
}
\DeclareAcronym{PI}{ 
    short = {PI}, 
    long  = {Proportional Integral Controller},
    tag = {abbrev}
}
\DeclareAcronym{FLC}{ 
    short = {FLC}, 
    long  = {Fuzzy Logic Controller},
    tag = {abbrev}
}
\DeclareAcronym{PID}{ 
    short = {PID}, 
    long  = {Proportional Integral Derivative Controller},
    tag = {abbrev}
}
\DeclareAcronym{UPFC}{ 
    short = {UPFC}, 
    long  = {Unified Power Flow Controller},
    tag = {abbrev}
}
\DeclareAcronym{THD}{ 
    short = {THD}, 
    long  = {Total Harmonic Distortion},
    tag = {abbrev}
}
\DeclareAcronym{PIDF}{ 
    short = {PIDF}, 
    long  = {PID Controller with Derivative Filter},
    tag = {abbrev}
}
\DeclareAcronym{ACO}{ 
    short = {ACO}, 
    long  = {Ant Colony Optimization},
    tag = {abbrev}
}
\DeclareAcronym{BA}{ 
    short = {BA}, 
    long  = {Bat Algorithm},
    tag = {abbrev}
}
\DeclareAcronym{SSA}{ 
    short = {SSA}, 
    long  = {Salp Swarm Algorithm},
    tag = {abbrev}
}
\DeclareAcronym{WOA}{ 
    short = {WOA}, 
    long  = {Whale Optimization Algorithm},
    tag = {abbrev}
}
\DeclareAcronym{FOPID}{ 
    short = {FOPID}, 
    long  = {Fractional Order Proportional Derivative Proportional Integral controller},
    tag = {abbrev}
}
\DeclareAcronym{SMC}{ 
    short = {SMC}, 
    long  = {Sliding Mode Controller},
    tag = {abbrev}
}
\DeclareAcronym{MPC}{ 
    short = {MPC}, 
    long  = {Model Predictive Control},
    tag = {abbrev}
}
\DeclareAcronym{RL}{ 
    short = {RL}, 
    long  = {Reinforcement Learning},
    tag = {abbrev}
}
\DeclareAcronym{PS}{ 
    short = {PS}, 
    long  = {Pattern Search algorithm},
    tag = {abbrev}
}
\DeclareAcronym{GA}{ 
    short = {GA}, 
    long  = {Genetic Algorithm},
    tag = {abbrev}
}
\DeclareAcronym{BH}{ 
    short = {BH}, 
    long  = {Black Hole algorithm},
    tag = {abbrev}
}
\DeclareAcronym{PSO}{ 
    short = {PSO}, 
    long  = {Particle Swarm Optimization},
    tag = {abbrev}
}
\DeclareAcronym{PIMR}{ 
    short = {PIMR}, 
    long  = {Proportional-Integral Multiresonant controller},
    tag = {abbrev}
}
\DeclareAcronym{MVO}{ 
    short = {MVO}, 
    long  = {Multi-Verse Optimizer},
    tag = {abbrev}
}
\DeclareAcronym{PID+DD}{ 
    short = {PID+DD}, 
    long  = {PID plus double-derivative controller},
    tag = {abbrev}
}
\DeclareAcronym{FPDPI}{ 
    short = {FPDPI}, 
    long  = {Fractional Order Proportional Derivative Proportional Integral controller},
    tag = {abbrev}
}
\DeclareAcronym{NN}{ 
    short = {NN}, 
    long  = {Neural Network},
    tag = {abbrev}
}
\DeclareAcronym{DE}{ 
    short = {DE}, 
    long  = {Differential Evolution},
    tag = {abbrev}
}
\DeclareAcronym{ABC}{ 
    short = {ABC}, 
    long  = {Artificial Bee Colony algorithm},
    tag = {abbrev}
}
\DeclareAcronym{GWO}{ 
    short = {GWO}, 
    long  = {Grey Wolf Optimizer},
    tag = {abbrev}
}
\DeclareAcronym{APF}{ 
    short = {APF}, 
    long  = {Artificial Potential Field},
    tag = {abbrev}
}
\DeclareAcronym{GSO}{ 
    short = {GSO}, 
    long  = {Galactic Swarm Optimization},
    tag = {abbrev}
}

\usepackage{multicol}

\AtBeginDocument{}

\hypersetup{colorlinks = true,linkcolor = blue,anchorcolor =red,citecolor = blue,
filecolor = red,urlcolor = red,pdfauthor=author}

\newcommand{\q}[1]{``#1''}

\journal{This is a preprint of a paper that is published open access 
in 'Heliyon' [\url{https://doi.org/10.1016/j.heliyon.2024.e31771}].}

\begin{document}

\begin{frontmatter}

\title{Universe-inspired algorithms for Control Engineering: A review}

\author[rodrigoAdress]{Rodrigo M. C. Bernardo\corref{correspondencia}}
\ead{rodrigo.bernardo@ua.pt}
\cortext[correspondencia]{Corresponding authors at: 
Department of Mechanical Engineering, University of Aveiro, 3810-193 Aveiro, Portugal 
(Rodrigo M. C. Bernardo and Marco P. Soares dos Santos).}

\author[delfimAdress]{Delfim F. M. Torres}
\ead{delfim@ua.pt}

\author[delfimAdress]{Carlos A. R. Herdeiro}
\ead{herdeiro@ua.pt}

\author[rodrigoAdress,lasi]{Marco P. Soares dos Santos\corref{correspondencia}}
\ead{marco.santos@ua.pt}

\address[rodrigoAdress]{Center for Mechanical Technology \& Automation (TEMA), 
Department of Mechanical Engineering,\\ 
University of Aveiro, 3810-193 Aveiro, Portugal}

\address[delfimAdress]{Center for Research and Development in Mathematics and Applications (CIDMA),
Department of Mathematics, \\
University of Aveiro, 3810-193 Aveiro, Portugal}

\address[lasi]{Intelligent Systems Associate Laboratory (LASI), Portugal}


\begin{abstract}
Control algorithms have been proposed based on knowledge related to nature-inspired mechanisms, including those based on the behavior of living beings. This paper presents a review focused on major breakthroughs carried out in the scope of applied control inspired by the gravitational attraction between bodies. A control approach focused on Artificial Potential Fields was identified, as well as
four optimization metaheuristics: Gravitational Search Algorithm, Black-Hole algorithm, Multi-Verse Optimizer, and Galactic Swarm Optimization. A thorough analysis of ninety-one relevant papers was carried out to highlight their performance and to identify the gravitational and attraction foundations, as well as the universe laws supporting them. Included are their standard formulations, as well as their improved, modified, hybrid, cascade, fuzzy, chaotic and adaptive versions. Moreover, this review also deeply delves into the impact of universe-inspired algorithms on control problems of dynamic systems, providing an extensive list of control-related applications, and their inherent advantages and limitations. Strong evidence suggests that gravitation-inspired and black-hole dynamic-driven algorithms can outperform other well-known algorithms in control engineering, even though they have not been designed according to realistic astrophysical phenomena and formulated according to astrophysics laws. Even so, they support future research directions towards the development of high-sophisticated control laws inspired by Newtonian/Einsteinian physics, such that effective control-astrophysics bridges can be established and applied in a wide range of applications. 
\end{abstract}


\end{frontmatter}


\section*{Abbreviations}
\printacronyms[include=abbrev, heading=None] 

\section{Introduction}

In recent decades, control of non-linear systems has been one of the most important topics in control theory \citep{Nasir2022}. Despite the massive use of non-linear models for accurate prediction of physical systems, it is still difficult to ensure high stability margins and desired performances in non-linear systems, mainly if uncertainties must be overcome \citep{Ali2018}. Researchers have been observing nature seeking inspiration to solve complex real-world control-related problems, since it is a clear example of a time-dependent process in a state of optimization, according to evolutionary mechanisms. One can find many natural processes in which a state of equilibrium and adaptation is reached, which can be investigated for nature-inspired high-performance optimization and control. Steer \textit{\textit{et al.}} \cite{steer2009rationale} stated that the term \textit{nature} refer "to any part of the physical universe which is not a product of intentional human design". These authors also distinguish between ‘strong’ inspiration and ‘weak’ inspiration, where the first one involves "the investigation of some existing problem-solving mechanism, the extraction of some qualitative process description, and the application to some alternative purpose", while the second is the "less formal role of some phenomenon in the creative stage of solution formulation".

Well known control methods do not consider the dynamics occurring in natural phenomena (non-inspired control) or only consider some dynamics occurring in biological structures (bio-inspired control). Many non-inspired control methods were already proposed, such as the Proportional, Integral and Derivative (PID) control, predictive control, optimal control, and sliding mode control \cite{fu2021hybrid,yu2020terminal,agarwal2018survey}. These are non-nature-inspired and employ "artificial" control approaches, often neglecting the rationality and effectiveness inherent in natural systems. In the case of sliding mode control, it does exhibit an attraction-like behavior, as the system state appears to be drawn towards the sliding surface. However, this attraction to the sliding surface is achieved through an artificial mechanism using a variable switching structure. Nonetheless, these controllers are formulated using a non-natural attraction, therefore they are not rooted in the natural behavior of celestial bodies in the universe. Intelligent control and bio-inspired control, including the design of Artificial Neural Networks and \ac{FLC} \cite{jin2020simulation,wu2018decentralized,yao2020model}, have also been extensively applied. However, their usage often demands a non-negligible degree of intuition and lacks interpretability \cite{ibrahim2022augmented,sarabakha2019intuit}. 

The main goal of this paper is to provide a literature review of the most relevant studies that highlight major scientific achievements in the domain of nature-inspired universe-conveyed control, to highlight their ability for future applications in multiple areas. Regarding their application in control systems engineering, scientific efforts have been centered on optimization and development of metaheuristics, despite the excess of metaphorical heuristics already reported \cite{aranha2022metaphor}. No control methods have been found with mathematical and physical formulations of gravitational attraction or black holes dynamics directly in their composition, thus evidencing a literature gap to be explored, where promising control methods may be designed using astrophysical dynamics, as they may provide mechanisms of stability and robustness (e.g. the strong gravitation field occurring in black holes). Moreover, the use of spacetime curvatures may hold great potential to engineer high performance trajectory tracking systems, as such phenomena ensure the shortest natural path between two points/states. Indeed, this review perform a thoroughly analysis to both optimization algorithms, already applied in the field of Control Engineering, and the control algorithms themselves, as long as their formulations are deeply related to gravitational attraction phenomena. This goal was achieved by providing an extensive analysis to the main concepts from which the original algorithms and related variants were developed, including their performance, characteristics and applications. 

After conducting an initial structured search, no actual control methods inspired by gravitational attraction were found. The closest approach involves the use of \ac{APF} to introduce attraction or repulsion behavior into systems. However, concerning optimization applied to control (\textit{e.g.} the optimization of controller parameters), four algorithms were identified: \ac{GSA}, \ac{BH}, \ac{MVO}, and \ac{GSO}. To our knowledge, no literature reviews were already focused on the use of gravitational phenomena in Control Engineering. Indeed, several review papers were already published in the scope of \ac{GSA}, \ac{BH}, \ac{MVO}, but they are mainly focused in data clustering, classification or general optimization problems \cite{rashedi2018comprehensive,abualigah2020multi}. Besides, bio-inspired control methods (including those inspired in swarm intelligence or evolutionary phenomena) are currently much more explored than non-biological nature inspired control methods, even though astrophysical phenomena hold potential to be used for developing high-sophisticated control systems, due to their inherent gravitational attraction.

The contributions of this paper can be summarized as follows: (i) Identification and exploration of control methods inspired by gravitational attraction or black-hole attraction derived dynamics; (ii) Critical analysis to the optimization algorithms already applied to control problems. Included are the \ac{GSA}, \ac{BH}, \ac{MVO}, and \ac{GSO}, as well as their variants and modifications; (iii) Critical analysis to the \ac{APF} already applied to control problems; (iv) Discussion on the potential advancements and limitations related to the use of gravitational attraction and universe-inspired algorithms in control systems. The ultimate goal is to contribute towards the development of high-sophisticated control systems inspired by realistic astrophysical phenomena and authentically formulated by Newtonian/Einsteinian physics.

\section{Methods}
\label{teoria}

\subsection{Selection Criteria}

In this paper we present a rigorous analysis of controllers and optimization algorithms applied in control systems inspired in a specific natural phenomenon: gravity, and related attraction between bodies. The Scopus database was searched in the time interval between 2000 and 2023 by seeking for the terms \q{gravit* AND control*}, \q{attrac* AND control*}, \q{black-hole AND control*}, and \q{galactic AND control*}, in the title, abstract and keywords. 
A search using the term \q{universe AND control} was also conducted; however, only studies outside the scope were obtained.
Control inspired in black-holes was included as they are currently considered an extreme phenomenon where extreme gravity and related extreme attraction conditions occur. The searching results were limited to: (i) document type: journals; (ii) subject: engineering; (iii) language: English. The compilation was further refined to remove documents outside the scope of this review, which as carried out according to the following rules:
\begin{enumerate}
    \item All the papers in control field obtained by searching the word "attraction", but referring to multiple meanings of the term not related to control science were removed (\textit{e.g.} interest, liking, and tempting).
    \item All the papers which contain the term "control" but do not refer to field of control systems (\textit{e.g.} attraction of ants by pheromones).
    \item All the papers that refer the terms "attraction", "gravity" and "black-holes" but whose controllers were not inspired in the gravity phenomenon (\textit{e.g.}, control in zero-gravity, micro-gravity environment).
    \item All the papers in the third and fourth quartiles, according to the Clarivate ranking, were removed, as we found they do not provide relevant content.
\end{enumerate}

The search was completed in March 2024. Ninety-one relevant papers were selected according to these criteria.

\subsection{Literature Search Strategy}
The following data were extracted and analyzed from the selected collection of papers: (1) inspired control law (concept, architecture, and analytical formulation); (2) inspired optimization algorithm (concept, type of optimization, and analytical formulation); (3) differences in the main concepts and main analytical formulations found in modified or hybrid versions when compared to the original proposed versions; (4) application of the proposed methods in the field of control systems; (5) relevant performance indicators.

\subsection{Terminology}

In the last decades, a large number of different nature-inspired algorithms and variants (\textit{e.g.} modifications and hybridization) were proposed to overcome relevant limitations mainly related to entrapment in local optima, premature convergence, parameter tuning, and exploration and exploitation imbalance \cite{thymianis2023hybridization}. Nevertheless, the adopted terminology to describe the different algorithms has not been widely consensual among researchers. 
Hence, concerning the conceptual differences between the original algorithm and its variants, the following classification was established:
\begin{itemize}
    \item \textbf{Standard:} The algorithm is used in its original formulation without any changes.
    \item \textbf{Improved:} The algorithm was upgraded  aiming to achieve superior performances, but without affecting the original conceptualization, (no artificial mechanisms were introduced).
    \item \textbf{Modified:} The algorithm was modified aiming to achieve superior performances,  but in such a way that partially or totally loses the affinity to its original conceptualization. 
    \item \textbf{Hybrid:} Merging of two algorithms aiming to achieve better performances in comparison with their individual performance. The algorithms must be truly combined in their formulation, i.e., they must not be formulated as an individual sequencing. 
    \item \textbf{Cascade:} Two algorithms individually sequenced. 
    \item \textbf{Fuzzy:} Algorithm that include  fuzzy logic in their conceptualization (\textit{e.g.},  for fine-tune parameterization). 
    \item \textbf{Chaotic:} Algorithm that include chaotic behaviors aiming to improve its performance.
    \item \textbf{Adaptive:}  Algorithm that include time-dependent modifications to the original algorithm throughout iterations (\textit{e.g.}, parametric modification).
\end{itemize}


\section{Attraction inspired optimization algorithms applied to control of dynamic systems}
\label{section:optimisation_alg}
\subsection{Gravitational Search algorithm}
\subsubsection{Overview}

\begin{figure}[ht]
\captionsetup{singlelinecheck = false, justification=justified}
\centering
\sidesubfloat[]{\label{subfig:figGSAintro}\includegraphics[width=0.45\linewidth]{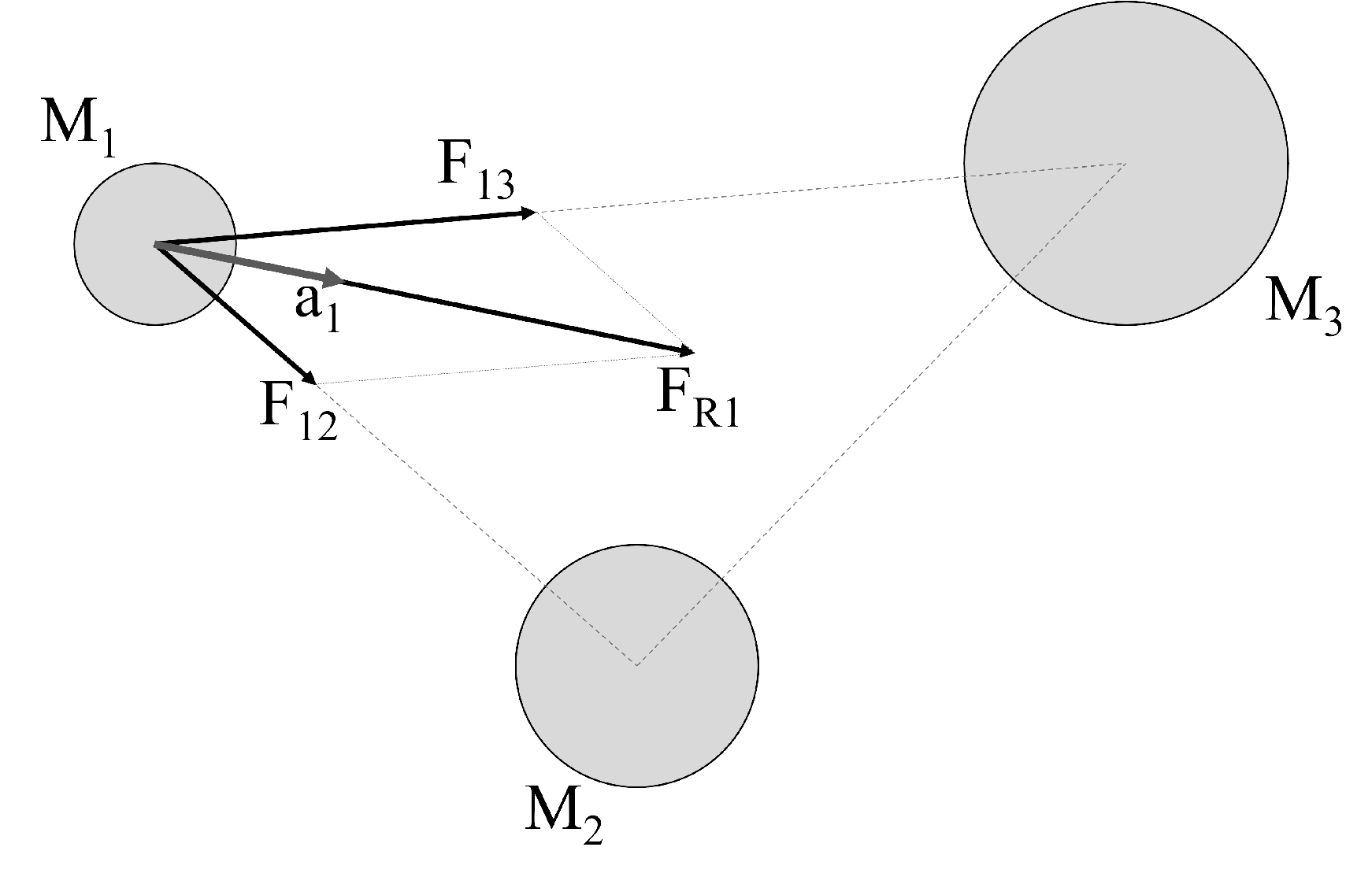}}
\hfill
\sidesubfloat[]{\label{subfig:figBHintro}\includegraphics[width=0.45\linewidth]{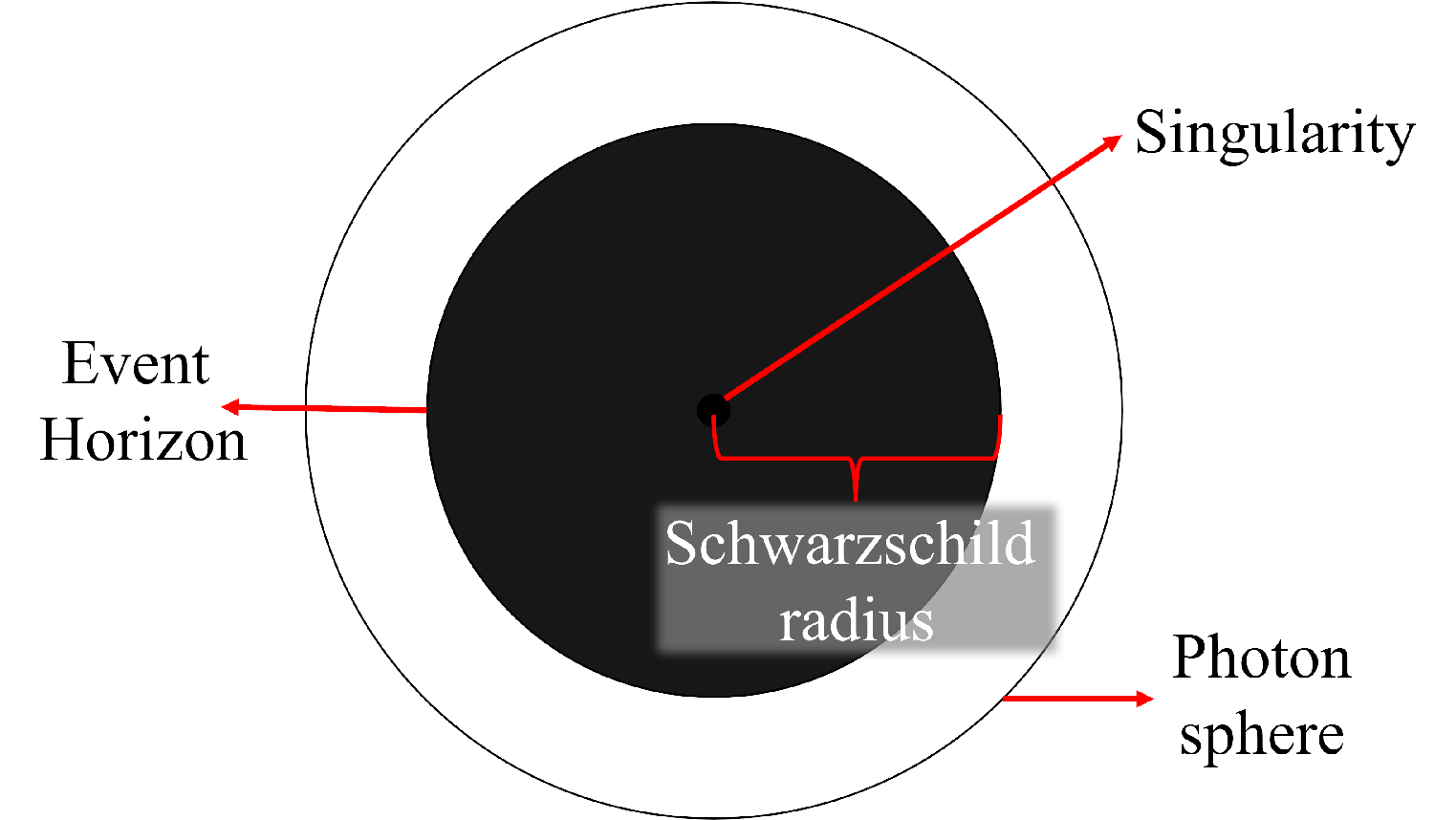}}  \\ [20pt]
\sidesubfloat[]{\label{subfig:MOV}\includegraphics[width=0.85\linewidth]{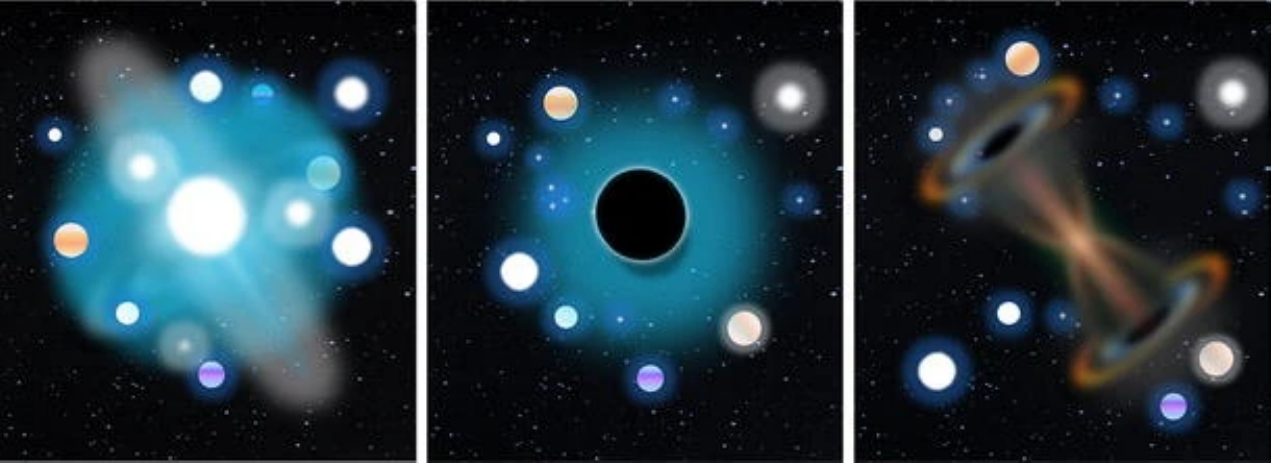}}
\caption{(a) Forces due to gravitational attraction on a three-body system. \(F_{12}\) 
is the force that \(M_{2}\) applies on \(M_{1}\), \(F_{13}\) is the force that 
\(M_{3}\) applies on \(M_{1}\), \(F_{R1}\) is the resultant force applied on \(M_{1}\) 
and \(a_{1}\) is the acceleration due to \(F_{R1}\); (b) Black Hole structure. 
The Schwarzschild \((R_{S})\) radius is calculated by \(R_{S}=\frac{2GM}{c^{2}}\) 
where \(M\) is the black hole mass, \(G\) is the gravitational constant and \(c\)
is the light speed. According to the Black Hole theories \cite{Kumar2015}, 
all objects that enter into the event horizon can not escape due to the 
massive gravitational attraction force. 
(c) Illustration of white-hole, black-hole and wormhole, respectively from left to right. Reproduced with permission from Ref.\cite{azizi2019optimal}
}
\label{fig:figura_introducao}
\end{figure}

The \ac{GSA} was firstly proposed by Rashedi \textit{\textit{et al.}} \cite{Rashedi2009}, who developed a heuristic optimization  method based on the Newton's law of gravity from classical physics. In this optimization algorithm, the search agents are represented by bodies whose mass depends on their fitness \cite{Rashedi2009}. The optimized solutions are obtained by body attraction phenomena, since bodies are modelled by larger masses to produce large attraction forces. Through this mechanism (\autoref{subfig:figGSAintro}) inspired by the gravitational force, the agents converge towards the best solution, which is represented by the body with the highest mass \cite{Rashedi2009}.
The baseline for the \ac{GSA} development was Newton's law of universal gravitation \cite{Newton1999},
\begin{align*}
\vec{\mathbf{F}}=-G\frac{M_1 M_2}{r^2} \hat{\mathbf{r}},
\end{align*}
where \(G\) is the gravitational constant, \(M_1\) and \(M_2\) 
are the bodies mass that attract each other, and \(r\) is the distance 
between the two bodies. Although the \ac{GSA} behaves as an artificial isolated system of masses with dynamics defined by the laws of gravitation and motion, these laws may be artificially modified from classic Newton law formulations, such that improved results can be achieved.

According to the original conceptualization \ac{GSA} \cite{Rashedi2009}, heavier bodies, which correspond to good solutions, move slowly than the lighter ones, which ensures the exploitation step of the algorithm. The method require to implement the \ac{GSA} as expressed in \autoref{subfig:figura_flowGSA}.

\captionsetup[figure]{justification = centering}
\captionsetup[subfigure]{justification = centering}

The optimization problem is modelled as a system with \(N\) mass-defined agents. The position of the \(i\)th agent is defined by
\begin{equation}
\label{eq_inicial}
X_{i} = (x_{i}^1,\ldots,x_{i}^d,\ldots,x_{i}^n) \;\;\;\;\; \textrm{for} \;i=1,2,\ldots,N
\end{equation}
where \(x_{i}^d\) is the position of the \(i\)th agent in the \(d\)th dimension. The gravitational constant \(G\) at time \(t\) is computed by
\begin{equation}
\label{eq_G}
G(t)= G_0 \textup{exp}\left ( -\beta\frac{t}{t_{max}} \right ) \; , \; \beta <1,
\end{equation}
where \(G_0\) is initial value of $G$. The Large Number hypothesis \citep{Dirac1974}, which was the first
hypothesis proposing a time varying gravitational constant, supported the paradigm stating that physical quantities should acquire dynamically their current values. Indeed, the \ac{GSA} was established by defining the force acting on mass \(i\) at time $(t)$ due to the presence of mass \(j\) as follows:
\begin{equation}
\label{eq:gsa_mainf}
F_{ij}^d(t)=G(t)\frac{M_{i}(t) \; M_{j}(t)}{R_{ij}(t)
+\varepsilon}\left( x_{j}^d(t)- x_{i}^d(t) \right)
\end{equation}
where \(\varepsilon\) is a small constant and \(R_{ij}(t)\) 
is the Euclidean distance between the two agents \(i\) and \(j\).
The total force acting on agent \(i\) in the dimension \(d\) is a
randomly weighted sum of \(d\)th components of the forces due to other agents:
\begin{equation}
\label{eq_totalForce}
F_{i}^d(t) = \sum_{j=1,j\neq i }^{N} rand_j F_{ij}^d(t),
\end{equation}
where \(rand_j\) is a random number in the interval $[0,1]$.
The position of the agents at the end of each iteration is calculated by:
\begin{equation}
\label{eq_pos}
x_i^d(t+1)=x_i^d(t)+v_i^d(t+1),
\end{equation}
\begin{equation}
\label{eq_vel}
v_i^d(t+1)=rand_i \; v_i^d(t)+a_i^d(t),
\end{equation}
\begin{equation*}
a_i^d(t) = \frac{F_i^d(t)}{M_i(t)},
\end{equation*}
where \(rand_i\) is a random value in the interval $[0,1]$; 
\(M_i(t)\) is the mass of \(i\)th agent at time \(t\), and it
is defined by
\begin{equation*}
M_i(t)=\frac{m_i(t)}{\sum_{j=1}^{N}m_j(t)},
\end{equation*}
\begin{equation*}
m_i(t)=\frac{fit_i(t)-worst(t)}{best(t)-worst(t)},
\end{equation*}
with \(fit_i(t)\) the fitness value of the \(i\)th agent at time \(t\), 
which depends on the defined objective function; \(best(t)\) and \(worst(t)\) 
are defined respectively by 
\begin{equation}
\label{eq_best}
best(t) = \left\{ 
\begin{array}{cl}
\min\limits_{j\in \left\{ 1,\dots,N \right\}} fit_j(t)  & \textrm{, if minimizing} \;\\
\max\limits_{j\in \left\{ 1,\dots,N \right\}} fit_j(t) & \textrm{, if maximizing}
\end{array} \right.
\end{equation}
and
\begin{equation}
\label{eq_worst}
worst(t) = \left\{ 
\begin{array}{cl}
\max\limits_{j\in \left\{ 1,\dots,N \right\}} fit_j(t)  & \textrm{, if minimizing} \;\\
\min\limits_{j\in \left\{ 1,\dots,N \right\}} fit_j(t) & \textrm{, if maximizing}
\end{array} \right.
\end{equation}
Functions (\ref{eq_best}) and (\ref{eq_worst}) are problem-dependent, 
i.e, minimization problems require a different formulation from maximization problems.

\subsubsection{GSA and related variations applied in control}

\begin{figure}[ht!]
\captionsetup{singlelinecheck = false, justification=justified}
\centering
\includegraphics[width=1\linewidth]{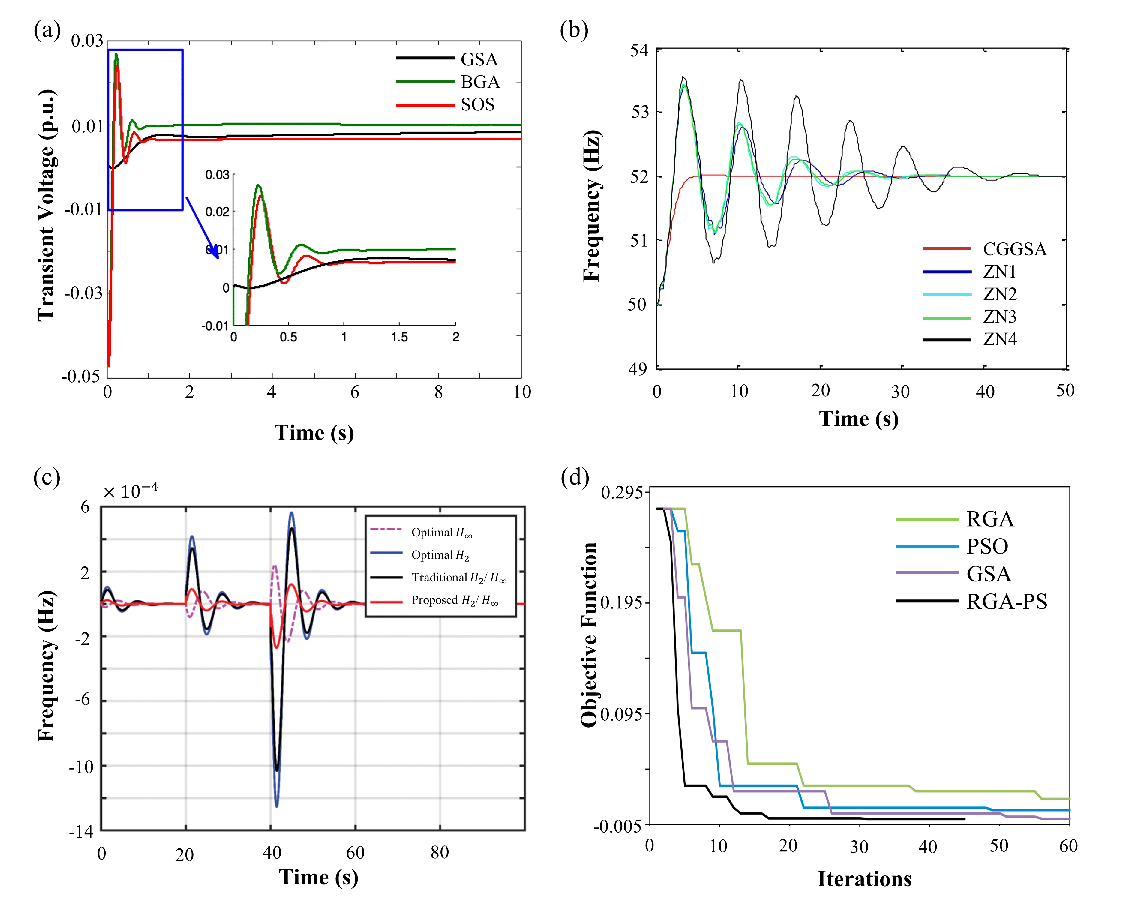}
\caption{(a) Control of transient voltage on a hybrid energy system using a PID controller tuned by GSA. Adapted with permission from Ref. \cite{guchhait2020stability}. (b) Frequency control of a pump turbine governing system using a \ac{PID} tuned by Ziegler–Nichols (ZN) method and by the proposed improved \ac{GSA} (CGGSA). Adapted with permission from Ref. \cite{li2017design}. (c) Comparison of different robust controller settings applied to control micro grid output frequency deviation, where the proposed $H_2/H_{inf}$ was optimized by hybrid particle swarm optimization and gravitational search algorithm with chaotic map algorithm (CPSOGSA). The proposed method was faster in retrieve the reference frequency with significantly less overshoot. Adapted with permission from Ref. \cite{zou2021optimized}. (d) Comparison of performance between Real Coded Genetic algorithm (RGA), PSO, GSA, and hybrid Real Coded Genetic - Pattern Search algorithm (RGA-PS). Adapted with permission from Ref. \cite{Chatterjee2014}.}
\label{fig:gsaimproved}
\end{figure}

Twenty-seven control applications were found  related to use of the standard version of \ac{GSA}, and thirty-six related to its variations (\autoref{table:originalGSA} and \autoref{table:variationsGSA}).
Applications of \ac{GSA} in control are mainly focused on optimal tuning of controllers gains, searching of the best control parameters, and finding the best control settings of complex systems. The main application field was electric energy generation (62\%), although they were already applied in the control of servo systems (10\%), as well as in applications with multiple constrains and requiring optimization of multiple parameters,  and also in applications in which control problems are transformed in optimization problems. The \ac{PID} (\autoref{fig:gsaimproved}a), \ac{FLC}, and \ac{UPFC}, whose parameters were optimized by some version of \ac{GSA}, represent the majority of the study cases (29\%, 16\%, and 5\%, respectively).

Concerning \ac{GSA} variants , the most used algorithm in control applications was the hybrid GSA-PSO, followed by chaotic mechanisms and improved versions of \ac{GSA}. Significant advantages have been found by using \ac{GSA} optimization algorithms. On the one hand, \ac{GSA} provides \citep{Khan2021,Sahu2015,Rashedi2011}: (i) a good global exploration capacity (good ability to search for new results); (ii) faster convergence in comparison to other methods (\textit{e.g.}, \ac{PSO}, \ac{GA} ; (iii) high computational efficiency; and (iv) higher accuracy in comparison to other methods (\textit{e.g.}, \ac{GA}, \ac{ACO} ,  On the other hand, limitations of \ac{GSA} are related to \citep{Zhao2018,Ji2017,Sun2018}: (a) diversity loss of new solutions in the final search steps;  (b) possibility of getting stuck in local optima; (c) parametrisation of the algorithm itself is required: their parameters have a significant influence in the effectiveness of the algorithm. Three modifications to \ac{GSA} were proposed so far to improve their effectiveness, by complementing the advantages of original \ac{GSA} with the advantages inherent to mechanisms of other searching or optimization methods: In order to find the optimal controller parameters of a hydraulic turbine governing system, an increasing $\beta$ value and a diversity based mutation were proposed \cite{li2016parameter}. The change performed on $\beta$ affects Eq. (\ref{eq_G}),  allowing to obtain a better control in the balance between exploration and exploitation. The second mechanism, triggered when the population diversity is lower than a dynamic threshold, ensures that the probability of agent mutation increases, such that the trap on local optima solutions is avoided. On the other two modified variants, changes were performed on Eq. (\ref{eq_vel}). To adjust the balance between global exploration and local exploitation, a simple mechanism based on a linear increasing $\gamma$ was introduced in Eq. (\ref{eq_vel}) to divide the equation in two terms. Concerning the problem of finding optimal \ac{UPFC} settings, an improvement of 2\% was achieved with less iterations in comparison with original \ac{GSA} \cite{deepa2017minimization}. Lu \textit{et al.} \cite{lu2019recurrent}  suggested a more complex  modification to the velocity update equation (\ref{eq_vel}) , including the transmission of information between agents to allow that all agents are updated based on the best ones, and adding memory to ensure that the best individual position is stored and used to compute (\ref{eq_vel}).   This concept is similar to the one used in \ac{PSO} \cite{kennedy1995particle}, despite it is differently formulated.

Some improved methods using non-complex concepts  were found to conduct to more effective results. To find the best thyristor controlled series compensator location to control a power system, Mahapatra \textit{et al.} \cite{mahapatra2016hybrid} proposed a mechanism to limit the maximum value of the velocity update (\ref{eq_vel}), with a decreasing maximum velocity, ensuring that the algorithm exploits the local search space in the final search phase. A similar approach was tested to optimize the thresholds and weights of a \ac{NN} model \cite{li2021intelligent}. Li \textit{et al.} \cite{li2017design} proposed to perform a mutation based on Gaussian and Cauchy distributions to enhance the exploitation and exploration capabilities of \ac{GSA}, respectively This method was tested by optimizing the controller gains of a pump turbine governing system, where the optimization capabilities of the improved \ac{GSA} were highlighted in contrast to the Ziegler-Nichols tuning approach (\autoref{fig:gsaimproved}b).

Opposition-based optimization is a technique already tested with many other optimization algorithms \cite{rojas2017survey,si2023pcobl,agarwal2021opposition,ewees2021new,shekhawat2020development,dhargupta2020selective}. Opposition optimization was used with \ac{GSA} and applied to control systems in order to find the optimal control parameters of power systems \cite{banerjee2014intelligent,shaw2014solution}. The main concept of opposition-based optimization is to check the opposite solution $\breve{x}_i$, defined as $\breve{x}_i = L + U - x_i$, where $L$ and $U$ are the lower and upper bounds of the search space, respectively. If the opposite candidate is fitter than the initial one, the opposite one is saved for the next iteration \cite{tizhoosh2005opposition}. Such optimization was also used with \ac{GSA}, and applied to control systems to find the optimal control parameters of power systems \cite{banerjee2014intelligent,shaw2014solution}.

Some processes, such as the \ac{GSA} tuning, are hard to determine objectively. However, Fuzzy Logic is a practical method of tuning the \ac{GSA} parameters as it can emulate the human reasoning in the use of imprecise information \cite{aghaie2017multi,precup2013fuzzy,diaz2017new}. Aghaie \textit{et al.} \cite{aghaie2017multi} proposed a fuzzy system to set the $\beta$ value in Eq. (\ref{eq_G}). Such proposed fuzzy system output new $\beta$ values according to four inputs: (1) the current iteration; (2) the progression level; (3) the diversification; and (4) the previous $\beta$ value. The diversification of population is given by

\begin{equation*}
    div = \frac{r_{ave} - r_{min}}{r_{max}-r_{min}}~,
\end{equation*}
where $r$ is the euclidean distance between two agents and $r_{ave}$, $r_{max}$ and $r_{min}$, are the average, maximum and minimum distances between agents, respectively. The level of progression is defined by
\begin{equation*}
    prog = \frac{fit_{ave}(t)-fit_{ave}(t-1)}{fit_{ave}(t)}~.
\end{equation*}
The proposed set of rules is shown in \autoref{table:fuzzyrulesAghaie}.

\begin{table}[ht]
\centering
\begin{tabular}{l l l l l l} 
 \hline
 \textbf{Rule} & \textbf{Inputs} & & & & \textbf{Output}\\ \cline{2-6}
    & Iteration & $prog$ & $div$ &$\beta(t-1)$ & $\beta(t)$ \\ \hline
 1 & low & low & low & medium & low \\
 2 & medium & low & low & high& medium \\
 3 & high & low & high & medium & high \\ \hline
\end{tabular}
\caption{Fuzzy rules used by Aghaie \textit{et al.} \cite{aghaie2017multi} to determine the $\beta$ values.}        
\label{table:fuzzyrulesAghaie}
\end{table}

Adaptation over iterations is other mechanism that has been employed by researchers to enhance the \ac{GSA} abilities. Two main conceptualizations were found using adaptive \ac{GSA} applied to control: by adapting $G$ and $\varepsilon$ values over time \cite{precup2012novel,precup2014adaptive}, and by performing a mutation with an adaptive probability, which is determined based on the success rate of the previous mutations \cite{niknam2013multiobjective}. Applications of deterministic chaos can be observed in control theory, computer science and physics; recently, chaotic-embedded \ac{GSA} has also been investigated  as another mechanism to improve the \ac{GSA} performance \cite{zelinka2022impact}. The use of chaotic maps allows to comprise additional layers of randomness to the algorithm, enhancing the local search capabilities \cite{li2013hydraulic,wei2018model}. By including neural behavior, Vikas and Parhi \cite{vikas2023multi} recently proposed a Modified Chaotic Neural Oscillator-based Hyperbolic \ac{GSA} (MHGSA) applied to humanoid robot path planning. They reported the ability of this adaptive GSA to achieve short paths in relation to original \ac{GSA} and avoid obstacles.

The most common hybrid algorithms applied to control is the hybrid GSA-PSO, due to the high similarities between \ac{GSA} and \ac{PSO} algorithms, which allows  an easy merging of their analytical formalization. Two versions were already proposed: (1) a simplified one only considering the propagation of the best solution through agents \cite{saha2017speed,mohanty2020novel,duman2015novel,khadanga2015time,deepa2018optimized,veerasamy2019automatic,khan2020design,ullah2020novel,vasant2020nature,veerasamy2022design}; and (2) a more complex one that saves the personal  best of each agent, adding memory to the algorithm \cite{bounar2019pso,kumar2022novel}. Other versions may arise through the combination of the various mechanisms mentioned above. Included is a chaotic hybrid GSA-PSO designed to optimize the parameters of a robust controller aiming to solve the load frequency problem of a micro grid, showing relevant results (improvement up to $\sim 83\%$), as shown in \autoref{fig:gsaimproved}c.

\begin{landscape}
\thispagestyle{empty}
\begin{longtable}{ >{\footnotesize\arraybackslash}p{3cm} >{\footnotesize\arraybackslash}p{16.5cm} c} 
\caption{Applications of original \ac{GSA} in control systems found in literature from 2000 to 2023. NA - Not  Applicable}        
\label{table:originalGSA}
\\
 \hline
 \textbf{Controller} & \textbf{Application Description} & \textbf{Reference}\\ [0.5ex] 
 \hline  
 NA & Find the optimal settings (\textit{e.g.} generator terminal voltages, transformer settings, output of compensating devices) for the reactive power dispatch problem that minimize the active power loss and enhance voltage stability of power system.& \cite{duman2012optimal} \\
 \acs{UPFC} & Search of optimal gains of \acs{UPFC} that exhibit greater robustness in the power system control. & \cite{ali2013multi} \\
 NA & Define the optimal switching angles of an inverter to minimize the \acs{THD}. & \cite{Chatterjee2014}\\
 State of charge feedback controller & Optimize the controller parameters to smooth the impact of photovoltaic sources in the power grid. &\cite{daud2014heuristic}\\
 \acs{PIDF} & Optimize the \acs{PIDF} gains of an Automatic Generation Control to minimize the generator frequency deviations and the tie-line power error of interconnected power systems. & \cite{Sahu2014} \\
 \acs{UPFC} & Search of optimal gains and location of multiple \acs{UPFC} that minimize the power loss and the dispatch cost of the power system. & \cite{sarker2014solution}\\ 
 \acs{FLC} & Find the optimal membership functions parameters of \acs{FLC}. The controller is applied to drive the speed of an induction motor. & \cite{abd2015novel}\\
  NA & Parametric optimization of ultrasonic machining processes& \cite{goswami2015parametric}\\
  \acs{PID}& Tuning of \acs{PID} gains to control a field-sensed magnetic suspension system & \cite{li2015gsa}  \\
  \acs{UPFC} & Find optimal settings of \acs{UPFC} during the post-fault period & \cite{kumar2015optimal} \\
  Type II/ Type III compensators & Find the optimal gains, zeros and poles location of the compensators to control a DC-DC boost converter & \cite{ghosh2016design}\\
  \acs{PI} & Optimize the \acs{PI} gains of an Automatic Generation Control to minimize the generator frequency deviations and the tie-line power error of interconnected power systems  & \cite{gupta2016performance} \\
  Fuzzy \acs{PID} & Optimize the controller parameters for Automatic Generation Control of a multi-area multi-source power system & \cite{pradhan2016firefly} \\
  \acs{FLC} & Optimize the rules and membership functions of \acs{FLC} to control the traffic flow & \cite{bi2017optimal} \\
  NA & Optimize the switching angles of a reactive power compensator & \cite{das2017reactive} \\
  Backstepping Control & Optimize the controller parameters for the trajectory tracking control of autonomous quadrotor helicopter & \cite{mohd2018trajectory}\\
  NA & Find the optimal settings to control the electric power generation system & \cite{ozyon2018gravitational} \\
  NA & Find the optimal electric vehicles controller settings that minimizes the voltage fluctuations and the degradation of batteries & \cite{ali2019voltage} \\
  \acs{PID} & Optimize the \acs{PID} gains to control an inverted pendulum system & \cite{Magdy2019} \\
  NA & Find the optimal settings of a congestion management system in a power system under deregulated regime & \cite{sharma2019gravitational}\\
  \acs{FOPID} & Optimize the \acs{FOPID} parameters to optimal control a micro grid system with various components & \cite{zaheeruddin2019design} \\
  \acs{SMC} & Optimize the \acs{SMC} parameters to control a dual-motor driving system & \cite{zeng2019fixed} \\
  \acs{PIDF} & Optimize the \acs{PIDF} gains to control a hybrid power system & \cite{guchhait2020stability} \\
  NA & Find optimal settings of a battery energy storage system & \cite{sakipour2020optimizing} \\
  \acs{MPC} & Optimize the \acs{MPC} parameters to determine online the optimal control sequence. Applied to a quadrotor & \cite{Nobahari2021} \\
  \acs{PI} & Determine the optimal parameters of a \acs{PI} to control the voltage and frequency of a micro grid & \cite{almani2022optimal} \\
  \acs{RL}-based control & Find optimal initial weights and biases of the Neural Network controller to avoid instability. The controller was tested in a linear position servo system & \cite{zamfirache2022reinforcement} \\
 \hline
\end{longtable}
\end{landscape}

\begin{landscape}
\pagestyle{empty}	
\begin{longtable}{ >{\footnotesize\arraybackslash}p{2cm} >{\footnotesize\arraybackslash}p{8.5cm} >{\footnotesize\arraybackslash}p{2cm} >{\footnotesize\arraybackslash}p{7cm} c} 
\caption{Applications of \ac{GSA} variations in control systems found in literature from 2000 to 2023. NA - Not  Applicable}        
\label{table:variationsGSA}
\\
 \hline
 \textbf{Variation} & \textbf{Modified Equations} & \textbf{Controller} & \textbf{Application Description} & \textbf{Reference}\\ [0.5ex] 
 \hline
 
 \multirow{9}{*}{\textbf{Modified}} & In (\ref{eq_G}) $\beta$ is updated as: $\beta(t) = \gamma \textup{sinh} \left( \eta\left( \frac{t}{t_{max}} - 0.5 \right) \right) + \lambda $ \newline 
 Perform the mutation: for each $x_i^d$, $x_i^d = z_i^d~\textup{if}~rand_1 < P_c ~\textup{where}~ Z_i = X_{r_1} + rand_2 (X_{r_2} - X_{r_3})$. The mutation occurs if $div < \varepsilon(t)~\textup{where}~\varepsilon(t) = \varepsilon_0 \textup{exp}\left( -\mu t / t_{max} \right)$ & \acs{PID} & Optimal parameter identification of a hydraulic turbine governing system & \cite{li2016parameter} \\ \cline{2-5}
 & Substitute (\ref{eq_vel}) by: $v_i^d(t+1) = \gamma \times gv_i^d(t+1)+(1-\gamma) \times lv_i^d(t+1)$ where $gv_i^d(t+1) = rand_i \times gv_i^d(t) + a_i^d(t)$ and $lv_i^d(t+1) = rand_i \times lv_i^d(t) + a_i^d(t)$& \acs{UPFC} & Find optimal \acs{UPFC} settings that minimize the power losses & \cite{deepa2017minimization}\\ \cline{2-5}
 & Substitute (\ref{eq_vel}) by: $v_i^d(t+1) = rand_0 v_i^d(t) + F \times a_i^d(t) + pbest_i^d(t) \left( 1-e^{-c_1 rand_1 \times t} \right) + gbest_i^d(t) \left( 1-e^{-c_2 rand_2 \times t} \right) $ where $F = \frac{1}{\sqrt{a}}e^{-\left( \phi / a \right)^2} \textup{cos} \left( 5\phi/ a \right)$&  Neural Network controller & Find neural network optimal parameters. The controller was applied on the integration of offshore wind and wave energy systems & \cite{lu2019recurrent} \\ \hline

 \multirow{12}{*}{\textbf{Improved}} & (\ref{eq:gsa_mainf}) is substituted by: \newline $F_{ij}^d(t)=G(t)\frac{M_{i}(t) \; M_{j}(t)}{R_{ij}(t)
+\varepsilon x_j^d(t)}\left( x_{j}^d(t)- x_{i}^d(t) \right)$ & \acs{FLC} & Optimize the \acs{FLC} parameters. The controller was applied on a DC servo system &  \cite{david2013gravitational} \\ \cline{2-5}
& Perform a velocity limitation in (\ref{eq_vel}): $-V_{max} \leq v_i^d \leq V_{max}$ where $V_{max} = V_{max0} \times \left [ 1-(t/t_{max})^h \right ]$ with $V_{max0} = \alpha (x_{max}-x_{min}), \alpha \in ]0,1]$ & NA  & Find the optimal thyristor controlled series compensator location in a power system & \cite{mahapatra2016hybrid} \\ \cline{2-5}
& Perform the following mutation to $X_i$: \newline $X_i^{new} = X_i \times (1 + \alpha \times (\eta N(0,1)+(1-\eta)C(0,1)))$ where $N(0,1)$ and $C(0,1)$ are random numbers from the Gaussian and Cauchy distributions respectively. Then, a new vector is obtained as $X_{all} = [X_{new} X]$ and only the best $N$ solutions are selected. & \acs{PID} & Optimize the \acs{PID} gains applied to a pump turbine governing system & \cite{li2017design} \\ \cline{2-5}
& Substitute (\ref{eq_vel}) by: $v_i^d(t+1) = rand_i w(t) v_i^d(t) +a_i^d(t)$, where $w(t) = w_{max} - \frac{w_{max}-w_{min}}{t_{max}}\times t$ & NA & Optimise the the neural network thresholds and weights. The neural network is used to filter the speed error that is used in the design of the servo system controller & \cite{li2021intelligent} \\ \hline

\multirow{5}{*}{\textbf{\shortstack[l]{Opposition\\based}}} & \multirow{5}{*}{\shortstack[l]{With a certain probability named jumping rate, $J_r$, after (\ref{eq_pos}) \\ the opposite solutions in relation to the actual population \\ are  verified: $Ox_i^d = \textup{min}^d + \textup{max}^d - x_i^d$. Then, the $N$ fittest\\ agents from set $\{ X, OX\}$ are selected.}} & \acs{FLC} & Search the optimal control parameters of an autonomous power system. The goal is to enhance the transient response, minimize the overshoot and oscillations, and improve the damping factor &  \cite{banerjee2014intelligent} \\ \cline{3-5}
 &  & NA & Search the optimal control parameters for the problem of optimal reactive power dispatch of power systems &  \cite{shaw2014solution} \\ \hline

\multirow{5}{*}{\textbf{\shortstack[l]{Fuzzy\\based}}} & The $\beta$ parameter in (\ref{eq_G}) is defined by fuzzy system. & NA & Search the optimal design parameters of core patterns for nuclear reactors to solve the loading pattern optimization problem &  \cite{aghaie2017multi} \\ \cline{2-5}
& The gravitational constant $G(t)$ and the parameter $\varepsilon$ in (\ref{eq:gsa_mainf}) are adapted using a fuzzy logic mechanism. & \acs{PI} & Search the optimal design parameters of a \acs{PI} controller for position control a servo system &  \cite{precup2013fuzzy} \\ \hline

\multirow{5}{*}{\textbf{Adaptive}} & \multirow{2}{*}{\begin{tabular}{l}
     On the first 15\% of iterations, (\ref{eq_G}) is given by: $G(t) = $\\
     $G_0 \left( \frac{1-\beta t}{t_{max}} \right)$. During next 45\% iterations (\ref{eq_G}) is defined by\\
     $G(t) = G_0 \textup{exp} \left( -\beta \frac{t}{t_{max}} \right)$ and $\varepsilon$ in (\ref{eq:gsa_mainf}) is given by $\varepsilon(t) =$\\
     $\varepsilon_0 - \varepsilon_0 \frac{(t - 0.15 t_{max})}{0.85 t_{max}}$. The remaining 40\% of iterations $G(t)$\\
     is set as a constant. 
\end{tabular}} & \acs{FLC} & Find optimal parameters of a \acs{FLC} for position control of a servo system \newline  & \cite{precup2012novel} \\ \cline{3-5}
 &  & \acs{PI} & Find optimal parameters of a \acs{PI} for position control of a servo system \newline   & \cite{precup2014adaptive} \\ \cline{2-5}
\multirow{5}{*}{\textbf{Adaptive}}& Two mutation mechanisms are considered: $X_{i,1} = X_{r_1} + rand_1 (X_{r_2}-X_{r_3}) + rand_2(X_{r_4}-X_{r_5})$ and $X_{i,2} = X_{r_1}+ rand_1 (X_{best} - X_{worst})$. On each iteration, the probability $P_a$ of occurring the mutation $X_{i,a}$ with $a=1, 2$ is given by $P_a = \frac{\textup{sr}_a}{\textup{sr}_1+\textup{sr}_2}$, where $\textup{sr}_a$ is the success rate of the mutation mechanism on past iterations. & NA & Applied for optimal reactive power dispatch and voltage control in power system operation & \cite{niknam2013multiobjective}  \\ \hline

\multirow{10}{*}{\textbf{Chaotic}} & \multirow{6}{*}{\shortstack[l]{After (\ref{eq_pos}), with a given probability, perform the chaotic search:\\
$ z_i(t) = \frac{x_{\textup{max}i} + x_{\textup{min}i}}{2} + \frac{x_{\textup{max}i} - x_{\textup{min}i}}{2} cx_i(t) $, where , $x_{\textup{max}i} =$\\
$x_{best} + \delta$ , $x_{\textup{min}i} = x_{best} - \delta$ and $cx_i(t)$ is a chaotic \\
map $\in [-1,1]\setminus 0$. If $z_i(t)$ is fitter than $x_{best}$ then $z_i(t)$ is the \\
new solution.}} & NA & Find the optimal parameters of a hydraulic turbine governing system fuzzy model \newline & \cite{li2013hydraulic} \\  \cline{3-5}
&  & Model-free controller & The optimal model-free controller was designed according to the quadratic performance index and applied to a vibro-impact system \newline &  \cite{wei2018model} \\ \cline{2-5}
&  \multirow{4}{*}{\shortstack[l]{(\ref{eq_totalForce}) is replaced by: $F_{i}^d(t) = \sum_{j=1,j\neq i }^{N} C(t) F_{ij}^d(t)$, \\
where $C(t)$ is a normalized chaotic map.}} & Robust controller & Optimal load frequency control settings applied to a micro grid &  \cite{zou2021optimized} \\  \cline{3-5}
&  & Robust \acs{FLC} & Find optimal controller parameters with application to the hydraulic turbine governing system &  \cite{li2022optimized} \\ \hline

\multirow{4}{*}{\textbf{\shortstack[l]{Chaotic neural\\ oscillators}}}& (\ref{eq_G}) is replaced by $G(t) = G_0 \times Z(t)$, with $Z(t) = 0.5(1 - L(t))$, where $L(t) = (V(t) - U(t))e^{-kt^2} + Y(t)$. Here, $Y(t) = \textup{tansig}(J(t))$. $U(t)$ and $V(t)$ are updated over iterations as following: $U(t+1) = \textup{tansig}\left( p_1 L(t) + p_2 U(t) - p_3 V(t) + p_4(J(t) - \phi_u) \right)$ and $V(t+1) = \textup{tansig}\left( q_1 L(t) - q_2 U(t) - q_3 V(t) + q_4(J(t) - \phi_v) \right)$. & NA & Find the optimal path for humanoid robot to avoid dynamic obstacles & \cite{vikas2023multi} \\ \hline

\multirow{5}{*}{\textbf{\shortstack[l]{Cascade}}}& After performing a global search using \ac{GSA}, the Gradient Descent Method is used to perform a refined local search (see \cite{ruder2016overview} additional details). & Lead-lag phase compensator & Find the optimal gains and lead-lag parameters of a power system stabilizer installed on synchronous generator & \cite{peres2018gradient}\\ \cline{2-5}
& After performing a global search using \ac{GSA}, the \acs{PS} is used to perform a refined local search (see \cite{sherif1994optimization} for additional details). & \acs{PID} & Find optimal automatic generation controller parameters to minimize frequency deviation of a multi-area electric power system & \cite{khadanga2017hybrid}\\ \hline

\multirow{2}{*}{\textbf{\shortstack[l]{Hybrid \\ GSA-FA}}}& (\ref{eq_pos}) is replaced by: $x_i^d(t+1) = x_i^d(t) + \alpha (x_j - x_i) + v_i^d(t+1)$, where $\alpha = \alpha_0 e^{- \gamma r^2}$ & \acs{PI} & Find optimal controller gains for the load frequency control problem of a power system & \cite{gupta2021hybrid}\\ \hline

\multirow{5}{*}{\textbf{\shortstack[l]{Hybrid \\ GSA-GA}}}& \multirow{5}{*}{\shortstack[l]{In each main iteration, $K\%$ of population is selected to evolve by\\
using \acs{GA} \cite{maier2019introductory}, and remaining population evolves using \acs{GSA}.\\
Then, solutions of both methods are combined. This process is\\
repeated until the maximum iteration is achieved.}} & \acs{UPFC} & Find optimal controller parameters to minimize system oscillations of a power system \newline & \cite{khadanga2015new}\\ \cline{3-5}
&  & \acs{FLC} & Find optimal controller parameters applied to speed control of a permanent magnet synchronous motor \newline & \cite{unsal2022investigation}\\ \hline

\multirow{6}{*}{\textbf{\shortstack[l]{Hybrid \\ GSA-PSO}}}& \multirow{6}{*}{\shortstack[l]{(\ref{eq_vel}) is replaced by:\\
$v_i^d(t+1) = w \times v_i^d(t) + k_1 rand_1 a_i^d(t)+ k_2 rand_2 (x_{best}^d - x_i^d(t))$}} & \acs{PI} & Find optimal speed controller parameters to minimize the ripple of a switched reluctance motor & \cite{saha2017speed}\\ \cline{3-5}
&  & \acs{FLC}/\acs{PID} & Find optimal automatic generation controller parameters to minimize the frequency deviation of a electric power system & \cite{mohanty2020novel}\\  \cline{3-5}
&  & NA & Find optimal settings of a multi-valve steam turbines system for power generation & \cite{duman2015novel}\\  \cline{3-5}
\multirow{19}{*}{\textbf{\shortstack[l]{Hybrid \\ GSA-PSO}}}& \multirow{13}{*}{\shortstack[l]{(\ref{eq_vel}) is replaced by:\\
$v_i^d(t+1) = w \times v_i^d(t) + k_1 rand_1 a_i^d(t)+ k_2 rand_2 (x_{best}^d - x_i^d(t))$}} & NA & Find optimal controller parameters of a power system stabilizer, applied to a multi-machine power system & \cite{khadanga2015time}\\  \cline{3-5}
&  & Neural Network \acs{MPC} & Find optimal neural network parameters of a non-linear continuous stirred tank reactor model used in the \acs{MPC} & \cite{deepa2018optimized}\\ \cline{3-5}
&  & \acs{PID} & Find optimal gains of \acs{PID} to control the interconnection of two area power system & \cite{veerasamy2019automatic}\\ \cline{3-5}
&  & NA & Find optimal settings for the optimal reactive power dispatch problem & \cite{khan2020design}\\ \cline{3-5}
&  & NA & The state estimation of a three-phase unbalanced distribution system is formulated as a nonlinear optimization problem which is solved by the proposed method & \cite{ullah2020novel}\\ \cline{3-5}
&  & NA & Applied to state-of-charge optimization (charging control) in the electric vehicles charging & \cite{vasant2020nature}\\ \cline{3-5}
&  & Neural \acs{PID} & Find optimal initial settings of the controller applied as automatic load frequency controller of interconnected hybrid power system & \cite{veerasamy2022design}\\ \cline{2-5}

& \multirow{6}{*}{\shortstack[l]{(\ref{eq_vel}) is replaced by:\\
$v_i^d(t+1) = w \times v_i^d(t) + k_1 rand_1 (x_{best,i}^d - x_i^d(t)) $\\
$+ k_2 rand_2 \left( \frac{a_i^d(t)+x_{best} }{2} - x_i^d(t) \right) $}} & Fuzzy \acs{SMC} & Find optimal controller parameters to control a generator-based wind turbine, ensuring power extraction maximization and regulation of reactive power according to grid requirements & \cite{bounar2019pso}\\ \cline{3-5}
&  & \acs{PID} & Optimal \acs{PID} for load frequency control of multi-source deregulated power system & \cite{kumar2022novel}\\

 \hline
\end{longtable}
\end{landscape}

\subsection{Black-Hole algorithm}
\subsubsection{Overview}

 A novel heuristic optimization method, motivated by the behavior of stars around a black hole (\autoref{subfig:figBHintro}), was proposed by Hatamlou \textit{\textit{et al.}}\cite{Hatamlou2013} in 2013: the \ac{BH}.  Its search agents are represented by stars; the one with the highest fitness value is established as the black hole-agent, which attracts all the others aiming to mimic the behavior of a real black hole. A star-agent is absorbed when it crosses the so-called Schwarzschild radius, which results in the remotion of the agent from the search space. To maintain a balance in the number of agents, a new star-agent is added at a random position of the search space. Throughout the iterations, if any star-agent becomes fitter than the black hole-agent, then the role of black hole-agent will be performed  in the next iteration.  After a predefined stopping criterion, the optimal solution is obtained by the the black hole-agent in the last iteration \cite{Hatamlou2013}. Even though the \ac{BH} is inspired by the behavior of the black hole phenomenon, it uses a conceptual form not supported by (Newtonian or Relativistic) physical laws  already theorized to describe the dynamics of black holes.

\begin{figure}[ht]
\captionsetup{singlelinecheck = false, justification=justified}
\centering
\sidesubfloat[]{\label{subfig:figura_flowGSA}\includegraphics[width=0.46\linewidth]{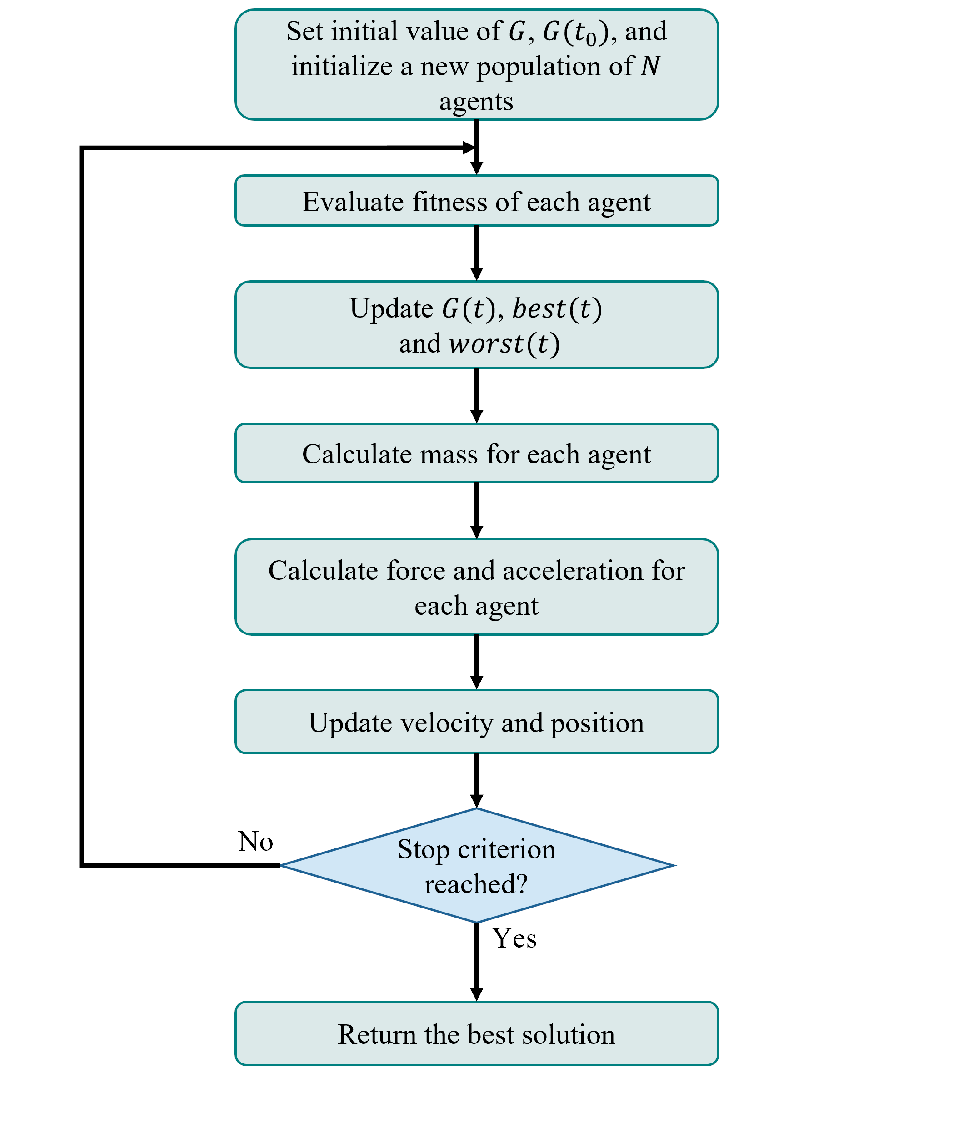}}
\hfill
\sidesubfloat[]{\label{subfig:figura_flowBH}\includegraphics[width=0.46\linewidth]{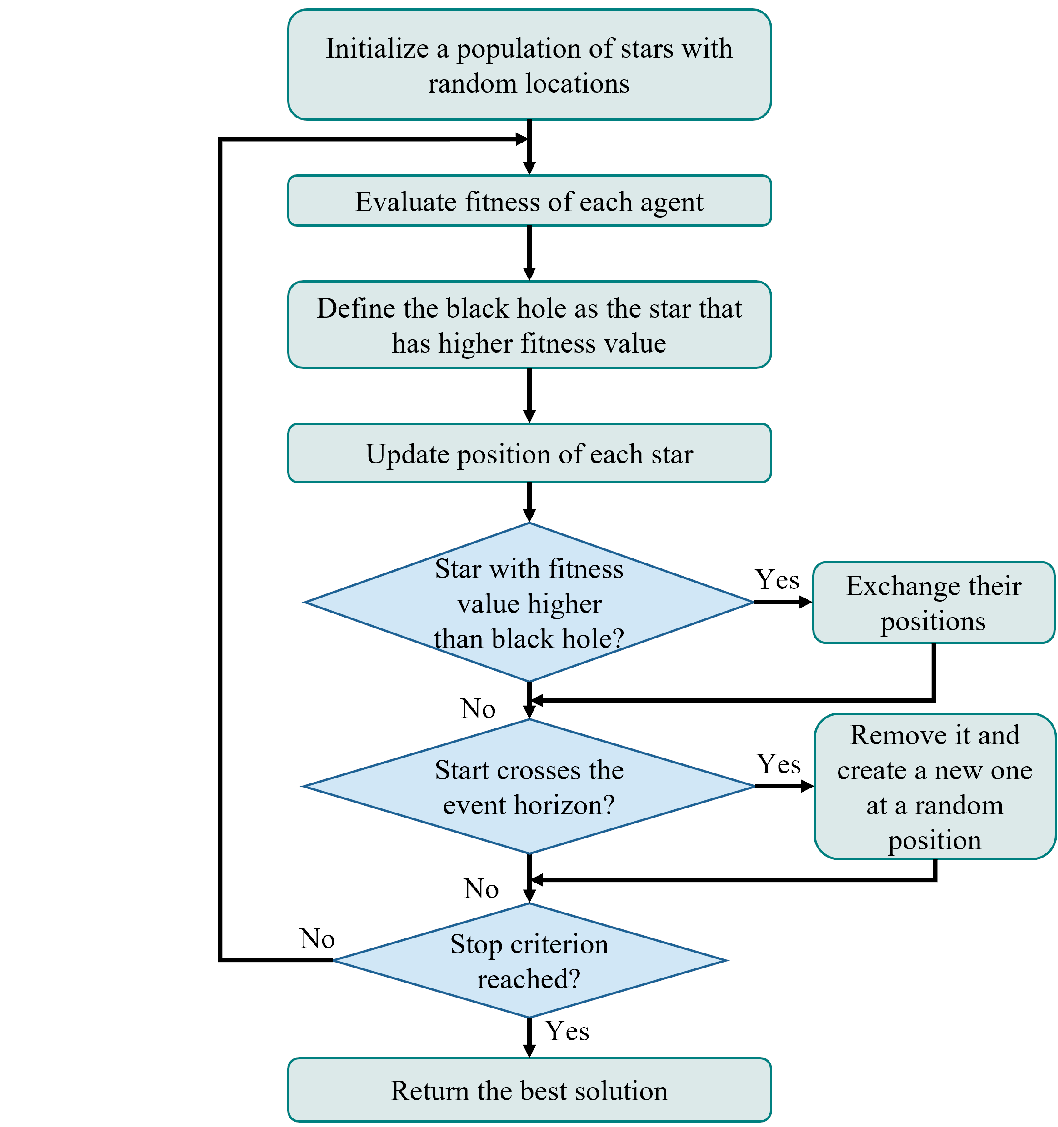}}\\
\sidesubfloat[]{\label{subfig:figura_flowMOV}\includegraphics[width=0.46\linewidth]{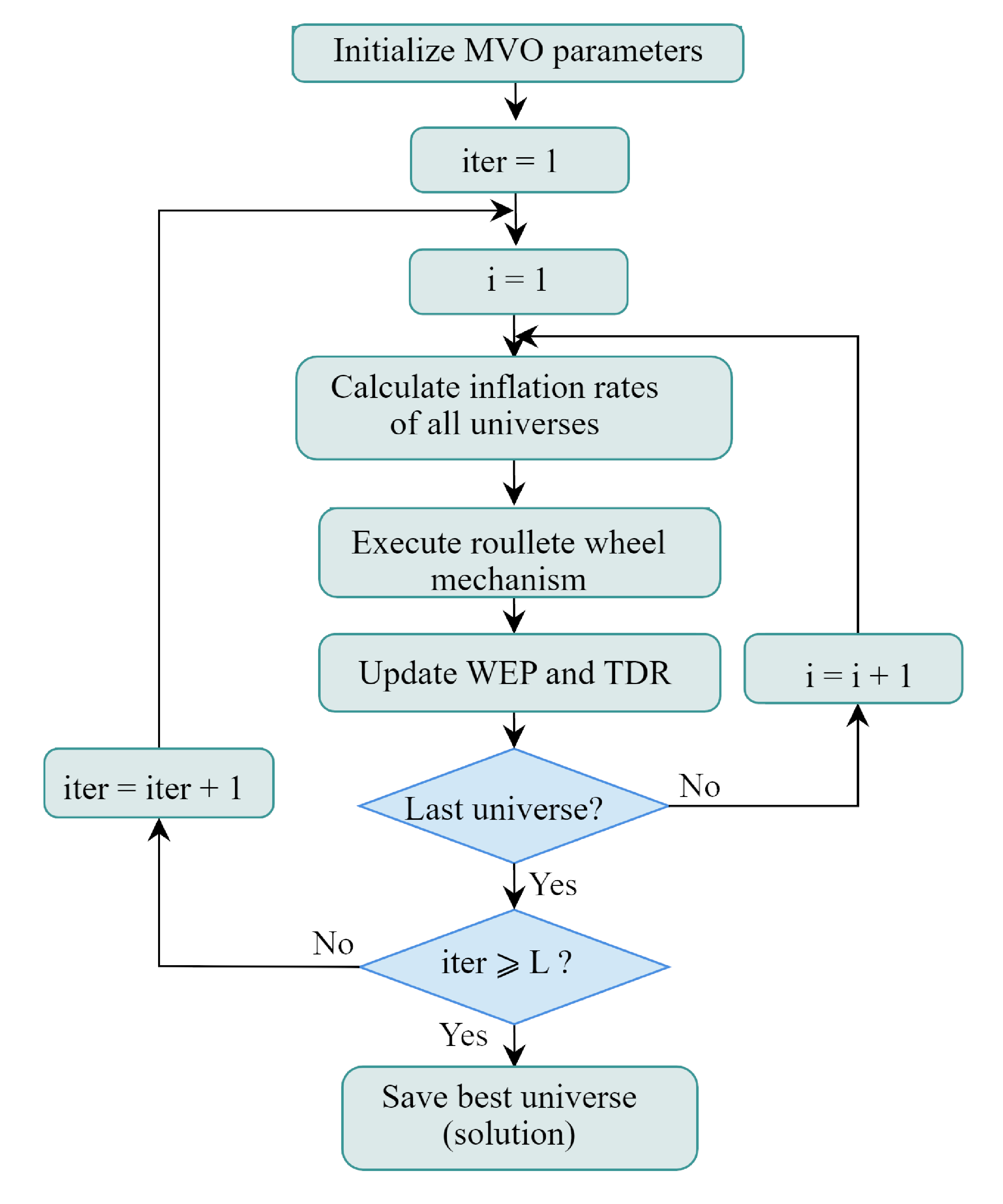}}
\hfill
\sidesubfloat[]{\label{subfig:figura_flowGSO}\includegraphics[width=0.46\linewidth]{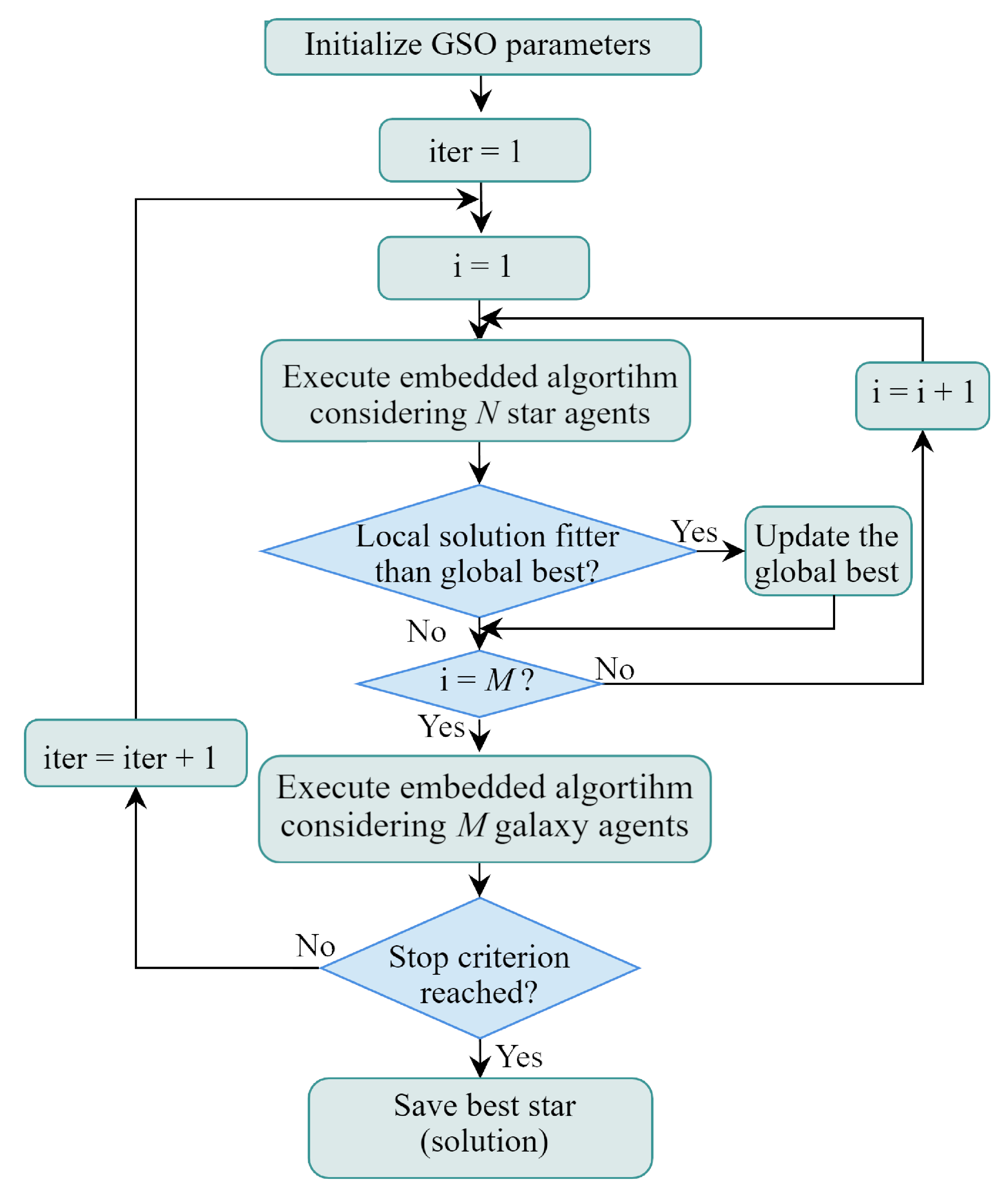}}
\caption{(a) Standard GSA algorithm. (b) Standard BH algorithm. (c) Standard MOV algorithm. (d) Standard GSO algorithm.}
\label{fig:flowcharts}
\end{figure}

Similarly to other population-based algorithms, the first step consists in generating an initial population of candidate solutions randomly distributed over the search space. Due to the small number of equations  formulating this algorithm, its implementation is not complex, even though its  high efficiency has been reported \cite{abualigah2022black}. The \ac{BH} algorithm is summarized in \autoref{subfig:figura_flowBH}. The original algorithm was established as follows. Let us consider a system with $N$ agents, in which the position of the $i$th agent is defined by Eq. (\ref{eq_inicial}). After the initialization, the fitness values of the agents are evaluated, and the best candidate is selected as the black hole-agent. As natural black holes absorb the stars surrounding them, the $i$th agent is dragged towards the black hole. The positioning of each agent is defined by

\begin{equation}
\label{eq_pos_bh}
x_i(t+1)=x_i(t)+rand \; \left( x_{BH} - x_i(t) \right),
\end{equation}
where \(rand\) is a random value in the interval $[0,1]$.
While mimicking the motion towards the black hole, if a star-agent becomes fitter than the black hole-agent in this new position, both this star-agent and the black hole-agent switch their positions. If during its movement, a star-agent crosses the event horizon of the black hole-agent, this star-agent \q{dies} and a new one is created randomly in the search space, to ensure a constant number of agents. The radius that defines the event horizon is given by
\begin{equation}
\label{eq_raio}
R = \frac{f_{BH}}{\sum_{i=1}^{N} f_i},
\end{equation}
where \(f_{BH}\) is the fitness value of the black hole 
and \(f_i\) is the fitness value of the \(i\)th star (agent).

\subsubsection{BH and variations applied in control}
Few applications using the \ac{BH} and its variations were reported, as  summarized in \autoref{table:BH}. The \ac{BH} was mostly applied in the optimization of parameters and gains of controllers (the exception was the modified BH). By using the original version of the \ac{BH}, the operation strategy (focused on finding the set point parameters) for a combined cooling, heating, and power system was optimized, such that the energy consumption, the system cost and the carbon dioxide emissions can be minimized \cite{deng2021novel}. Only a comparison with a \ac{PSO} algorithm was conducted: while the \ac{PSO} algorithm achieved optimal results with an objective function value minimized up to 0.595, the \ac{BH} obtained was minimized up to 0.58, which represents the  slight improvement of 2.5\%. The original algorithm was also applied to enhance the power quality of an AC micro grid by searching the optimal \ac{PIMR} gains that minimize the Total Harmonic Distortion (\ac{THD}) \cite{khosravi2022new}. The comparison was also carried out with \ac{PSO}, but higher performances were observed, namely 33\% improved objective function value and a faster convergence. 

An improved version of the \ac{BH} was found aiming to improve the motion of agents by introducing a new concept that prevents the dispersion of solutions: instead of a randomly generation, the obtained data of existing agents is used to generate new members(\autoref{table:BH}) \cite{azizipanah2016multiobjective}. This version was able to achieve good results when used to optimize the parameters of membership functions of a \ac{FLC}. By computing the minimization of the carbon emissions, this method was able to provide improvements of 17\% and 14\%, over the \ac{GA} and \ac{PSO}, respectively \cite{azizipanah2016multiobjective}. A similar formulation was used for urban traffic network control,  leading to improvements of 29\% compared to an already well-established approach \cite{Khooban2017}. This improved version was also applied to optimize the parameters of a Model-free \ac{SMC} for a Frequency Load Controller, designed to regulate a micro grid \cite{khooban2017secondary}. The proposed approach provided the best regulation under load changes.

\captionsetup[figure]{justification = centering}
\captionsetup[subfigure]{justification = centering}

\begin{figure}[ht]
\captionsetup{singlelinecheck = false, justification=justified}
\centering
\sidesubfloat[]{\label{subfig:xie2020}\includegraphics[width=0.85\linewidth]{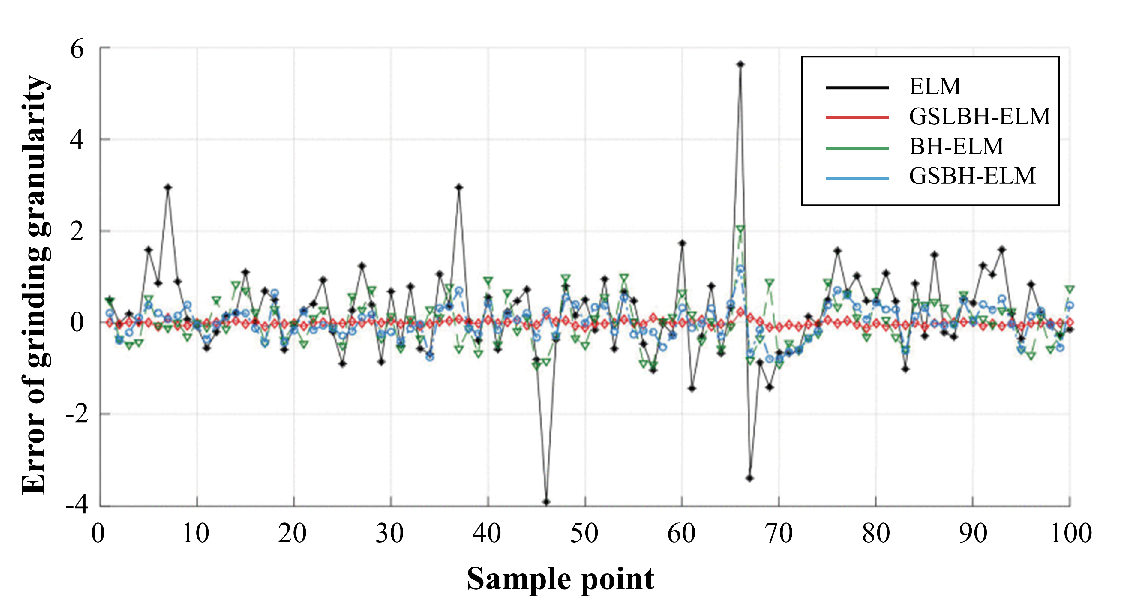}}
\\
[40pt]
\sidesubfloat[]{\label{subfig:aziziOptimal}\includegraphics[width=0.95\linewidth]{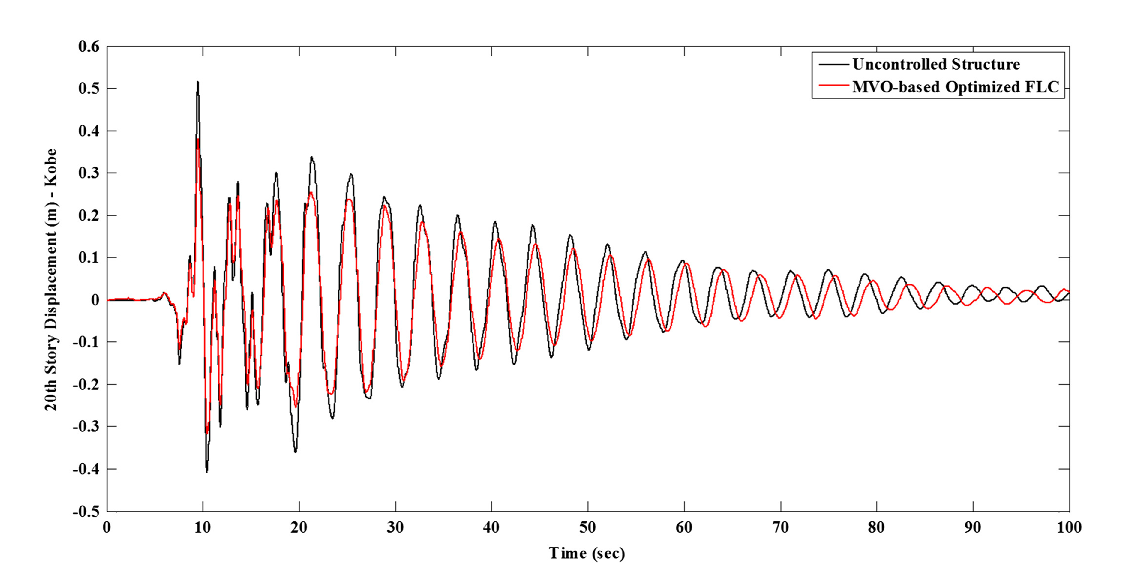}}
\caption{(a) Error of grinding granularity using the proposed soft-sensor model with different optimization algorithms. Note that the Golden Sine Levy-flight BH (GSLBH-ELM) achieved better results than without Levy-flight operator (GSBH-ELM) and than the base algorithm (BH-ELM). Adapted with permission from Ref. \cite{xie2020extreme}. (b) Displacement of a structure during an earthquake with and without structure control (MVO-based Optimized FLC). Reproduced with permission from Ref. \cite{azizi2019optimal}.}
\label{fig:bhversions}
\end{figure}

A modified version was proposed to optimize an extreme learning machine soft-sensor model to predict the grinding granularity \cite{xie2020extreme}. Comparing to the original \ac{BH} algorithm, this modified \ac{BH} was upgraded by applying two well-known operators to the movement of the agents, namely the Golden Sine operator \cite{li2024improved, li2023local,han2022golden,zhang2020improved}, and the Levy flight operator \cite{chawla2018levy,ewees2022improved,liu2021improving}. These operators have already been used to modify other optimization algorithms, namely \ac{BA}, \ac{SSA}, \ac{GWO} and \ac{WOA} \cite{li2023local,zhang2020improved,xie2013novel,zhang2020improved2}, and its ability to improve the original \ac{BH} algorithm has been recently demonstrated. Lower prediction errors were achieved in comparison to other methods, namely the original \ac{BH} and the Golden Sine \ac{BH} (without the Levy flight operator) (\autoref{subfig:xie2020}).

\begin{landscape}
\thispagestyle{empty}
\begin{longtable}{ >{\arraybackslash}p{2cm} >{\arraybackslash}p{8cm} >{\arraybackslash}p{3cm} >{\arraybackslash}p{6.5cm} c} 
\caption{Applications of \ac{BH} and variations of it in control systems found in literature from 2000 to 2023. NA - Not Applicable}        
\label{table:BH}
\\
 \hline
 \textbf{Variation} & \textbf{Modified Equations} & \textbf{Controller} & \textbf{Application Description} & \textbf{Reference}\\ [0.5ex] 
 \hline
 \rule{0pt}{3ex}   
\multirow{6}{*}{\footnotesize \textbf{Standard}}  & & & & \\
 & \footnotesize \centering NA & \footnotesize NA & \footnotesize Optimization of the set points of a combined cooling, heating, and power-ground source heat pump system. & \cite{deng2021novel}\\ \cline{3-5}
 & \footnotesize \centering NA & \footnotesize \acs{PIMR}  &\footnotesize Optimization of the controller gains. The controller was applied to improve the power quality components of an AC microgrid.  & \cite{khosravi2022new}\\
 \hline

 \multirow{6}{*}{\footnotesize \textbf{Modified}} & & & & \\
 & \footnotesize After (\ref{eq_pos_bh}) perform the Golden Sine and Levy flight operators: \newline $x_i(t+1) = x_i(t) + a\cdot \textup{sign}(rand - 0.5)\otimes s$ \newline $x_i(t+1) = \left | \textup{sin}(r_1) \right |x_i(t) - r_2\textup{sin}(r_1)\left | m_1 x_{BH}-m_2 x_i(t) \right |$ where $a$, $m_1$ and $m_2$ are coefficients, $s$ is the step size vector for Levy flight, $r_1$ is a random number in $[0,2\pi]$ and $r_2$ is a random number in $[0,\pi]$ & \footnotesize NA & \footnotesize Optimization of extreme learning machine (ELM) soft-sensor model on grinding process & \cite{xie2020extreme} \\
 \hline 
 
 \multirow{10}{*}{\footnotesize \textbf{Improved}} & & & & \\
  & \multirow{8}{*}{\shortstack[l]{\footnotesize (\ref{eq_pos_bh}) is changed for: $x_i(t+1) = x_i(t)+rand_1(x_{BH}-x_i(t))$ \\ \footnotesize $+rand_2(x_r(t)-x_i(t))$ where $r \in [1,N]$ \\
 \footnotesize (\ref{eq_raio}) is changed for: $R = \frac{\sum_{i=1}^{N} \left \| x_i - x_{mean} \right \|}{N} $ \\
 \footnotesize New candidate solutions are created by: $x_{new} = x_{BH} + $ \\ \footnotesize $(2rand(1,N)-1)\frac{\textup{max}(\left \| x_i - x_{BH} \right \|)}{N}$}} & \footnotesize \acs{FLC} &  \footnotesize Applied to \acs{FLC} membership functions parameters optimization. The controller was applied in a multi-objective dynamic optimal power flow framework. & \cite{azizipanah2016multiobjective} \\ \cline{3-5}
   &  & \footnotesize Model-Free \acs{SMC} &  \footnotesize Applied to optimize the parameters of Model-Free \acs{SMC}. The controller is then applied as Frequency Load Controller in a microgrid. & \cite{khooban2017secondary} \\ \cline{3-5}
   &  & \footnotesize \acs{FLC} &  \footnotesize Applied in control of traffic signal scheduling and phase succession to ensure smooth traffic flow with the objective of minimize the waiting time and average queue length. & \cite{Khooban2017} \\
[5ex] 
 \hline
\end{longtable}
\end{landscape}

\subsection{Multiverse algorithm}
\subsubsection{Overview}
The \ac{MVO} is a recent population-based optimization algorithm inspired in the multi-verse theory, focused on the interaction between universes, from which  white-holes, black-holes and wormholes emerge \cite{mirjalili2016multi}. While white-holes present similarities with  universes under expansion, black-holes attract everything with their extreme gravitational force. Wormholes are responsible for connecting different parts of a universe, acting like space-time traveling tunnels. These three cosmic objects of \ac{MVO} are illustrated in \autoref{subfig:MOV}.

 Each solution of the \ac{MVO} is analogous to a “universe”, and each solution dimension is an object that can be transmitted through “white-holes”, “black-holes” and “wormholes”. The objects are transferred from “white-holes” of a source-universe-solution to “black-holes” of a destination-universe-solution. Therefore, the population, corresponding to the set of universes-solutions, is described as
\begin{equation*}
    U = \begin{bmatrix}
x_1^1 & x_1^2 & \cdots   & x_1^d \\ 
x_2^1 & x_2^2  & \cdots & x_2^d \\ 
\vdots  &\vdots  & \ddots  & \vdots \\ 
 x_n^1 & x_n^2 & \cdots & x_n^d
\end{bmatrix}
\end{equation*}
where $d$ is the dimension of search space, and $n$ is the number of candidate solutions. For each object $x_i^j$, which denotes the $j$th variable of $i$th universe-solution, the following comparison is performed:
\begin{equation*}
    x_i^j = \left\{\begin{matrix}
x_k^j & r_1 < NI(U_i) \\ 
x_i^j & r_1 \geq  NI(U_i)
\end{matrix}\right.
\end{equation*}
where $NI(U_i)$ is the normalized inflation rate of the $i$th universe, $r_1$ is a random number in $[0,1]$, and $k$ is a universe-solution selected by a roulette wheel selection mechanism \cite{mirjalili2016multi}. The inflation rate of a universe-solution is a value proportional to the fitness of the corresponding solution.
This mechanism performs the exchange of objects between universes-solutions; in order to provide local changes, wormhole tunnels are established between a specific universe-solution and the best universe-solution emerged at time $t$. The formulation of such mechanism is the one that follows:
\begin{equation}
\label{eq:r4mov}
    x_i^j = \left\{\begin{matrix}
\left\{\begin{matrix}
X_j+ \textup{TDR}\cdot (r_4\cdot (ub_j-lb_j)+lb_j) & r_3 <0.5 \\ 
X_j- \textup{TDR}\cdot (r_4\cdot (ub_j-lb_j)+lb_j) & r_3 \geq 0.5
\end{matrix}\right. & r_2< \textup{WEP}\\ 
x_i^j & r_2\geq \textup{WEP}
\end{matrix}\right.
\end{equation}
where $X_j$ is the $j$th variable of the best solution, TDR (traveling distance rate) and WEP (wormhole existence probability) are coefficients, $lb_j$ and $ub_j$ are respectively the lower and upper bounds of the $j$th variable, and $r_2$, $r_3$ and $r_4$ are random numbers defined in $[0,1]$. The \ac{MVO} algorithm is illustrated in \autoref{subfig:figura_flowMOV}.

\subsubsection{MVO applied in control}

Only this standard \ac{MVO} was proposed so far, and its application is reduced to optimization of controller parameters, namely PID-derived and \ac{FLC} controllers (\autoref{table:MOV}). The first application in control systems of this method dates back to 2017, which was engineered to search for the optimal parameters of a \ac{PID+DD} to operate  as a Load Frequency Controller on a power system \cite{guha2017multi}. Recently, the same problem was revisited using the same algorithm to optimize a \ac{FPDPI}. \ac{MVO} was also applied in the active structural control of a building (vibration control) by searching the optimal set of parameters for the membership functions of a \ac{FLC} \cite{azizi2019optimal}, which reduced the vibration of a structure during tests with Kobe earthquake data, as shown in \autoref{subfig:aziziOptimal}. The performance of the \ac{MVO} was compared with other optimization metaheuristics (namely \ac{GA}, \ac{PSO} and \ac{GWO}) using test functions. However, in all the aforementioned papers, no comparisons were found extending to other optimization methods considering real case problems. 

\begin{table}[ht!]
\centering
\begin{tabular}{c >{\arraybackslash}p{4cm} >{\arraybackslash}p{7.5cm} c} 
 \hline
 \textbf{Variation} & \textbf{Controller} & \textbf{Application Description} & \textbf{Reference}\\ [0.5ex] 
 \hline
 \rule{0pt}{3ex}   
\multirow{6}{*}{Standard} & \footnotesize \acs{PID+DD} & \footnotesize Optimization of controller parameters. The controller was applied as a load frequency controller used to control the flow of steam to the turbines of a generator. & \cite{guha2017multi}\\ \cline{2-4}
 & \footnotesize \ac{FLC} &\footnotesize Optimization of the membership functions. The FLC was applied in active control of structures in civil engineering.  & \cite{azizi2019optimal}\\ \cline{2-4}
 & \footnotesize \acs{FPDPI} & \footnotesize Optimization of the controller parameters. Controller applied in a multi area power system consisting in hydro, thermal, and gas power plants. &\cite{sahoo2022multi} \\
[5ex] 
 \hline
\end{tabular}
\caption{Applications of MOV in control systems found in literature from 2000 to 2023.}        
\label{table:MOV}
\end{table}

\subsection{Galactic Swarm Optimization algorithm}
\subsubsection{Overview}
To enhance the equilibrium between exploration and exploitation, the \ac{GSO} was introduced in 2016 by Muthiah-Nakarajan and Noel \cite{muthiah2016galactic}. This algorithm draws inspiration from the movement of galaxies and the stars within them. Stars are not uniformly distributed throughout the cosmos; rather, they cluster into galaxies, which are not evenly distributed.

While \ac{GSO} has been conceptualized based on \ac{PSO}, the authors emphasize that this choice was made primarily due to the simplicity of implementing \ac{PSO}. They assert that \ac{GSO} could be implemented using any population optimization heuristic \cite{muthiah2016galactic}. Hence, the base method is not delineated to maintain generality.

To implement \ac{GSO}, $M$ galaxies are created, each containing $N$ different stars. During each iteration of the algorithm, the core heuristic is executed for each galaxy to determine the optimal solution within each galaxy. If a superior solution to the current global one is discovered during this process, the global solution is updated. At the conclusion of this phase, each galaxy is represented by its best local solution. In the subsequent phase, the core heuristic is applied, with the galaxies acting as the search agents, moving towards the best solution. This iterative process continues until a stopping criterion is met.

In summary, \ac{GSO} entails the application of a fundamental algorithm to ascertain the best local solution for each galaxy, followed by applying the algorithm at a broader level to determine the best global solution. Analogously, it can be conceptualized that during the initial phase, stars converge within each galaxy toward the star with the greatest mass, whereas in the subsequent phase, galaxies (clusters of stars) converge toward the galaxy with the greatest mass.
The \ac{GSO} algorithm is illustrated in \autoref{subfig:figura_flowGSO}.

\subsubsection{GSO applied in control}

The utilization of \ac{GSO} in system control remains relatively limited and under explored: only five distinct applications of this algorithm were identified. Among these cases, only one notably study employed the \ac{WOA} algorithm as its foundation, while the remaining cases utilized \ac{PSO}, as outlined in \cite{muthiah2016galactic}. These applications encompass the optimization of \ac{FLC} membership functions \cite{bernal2019optimization,bernal2021optimization}, the refinement of micro grid parameters \cite{gajula2020agile,rajasekaran2021bidirectional}, and the determination of optimal control parameters for an IoT network \cite{karthick2022galactic}. A summary of the research findings is presented in \autoref{table:GSO}.

\begin{table}[ht!]
\centering
\begin{tabular}{c >{\arraybackslash}p{3cm} >{\arraybackslash}p{7.5cm} c} 
 \hline
 \textbf{Embedded algorithm} & \textbf{Controller} & \textbf{Application Description} & \textbf{Reference}\\ [0.5ex] 
 \hline
 \rule{0pt}{3ex}   
\footnotesize PSO & \footnotesize \acs{FLC} & \footnotesize Optimization of the membership functions of the controller. The \ac{FLC} was applied to an autonomous mobile robot for trajectory tracking under noise effects. & \cite{bernal2019optimization,bernal2021optimization}\\[2ex] \cline{1-4}
\footnotesize PSO & \footnotesize \ac{FLC} &\footnotesize Optimization of the membership functions. The FLC was applied to control the water level in a water tank.  & \cite{bernal2020fuzzy}\\ \cline{1-4}
\footnotesize PSO & \footnotesize NA & \footnotesize Optimization of the settings parameters to provide the optimal power flow of a hybrid energy micro grid. &\cite{gajula2020agile} \\
\cline{1-4}
\footnotesize PSO & \footnotesize PI & \footnotesize Optimization of a Maximum Power Point Tracking controller to control a micro grid composed by photovoltaic panels and wind generators.&\cite{rajasekaran2021bidirectional}  \\[2ex] 
\cline{1-4}
\footnotesize PSO-WOA & \footnotesize NA & \footnotesize Finding the optimal settings for IoT network. &\cite{karthick2022galactic} \\[2ex] 
 \hline 
\end{tabular}
\caption{Applications of GSO in control systems found in literature from 2000 to 2023. NA - Not Applicable}        
\label{table:GSO}
\end{table}

\section{Attraction Phenomena Applied to Control of Dynamic Systems}

\subsection{Overview}

Attraction is a fundamental phenomenon that governs a wide range of interactions across the universe. In both physical and engineering senses, attraction and gravitational potential embody the notion of entities drawn towards one another, whether in the physical space or within the complex relationships of intelligent natural entities.
In the context of swarm intelligent systems, attraction plays a significant role in coordination of social organisms, including ants, bees, and birds, embedding the principles of collaboration, emergence, and decentralized decision-making. On the cosmos, gravity is the primary phenomena responsible for the attraction between bodies, formation of stars and planets, as well as to maintain stable orbits.

Inspired in the attraction phenomena, in the last decade, techniques to control swarms of intelligent systems have been developed, aiming to control both the swarm dynamics and its formation \cite{zhang2018fixed,liu2018distributed}. The proposed controllers were designed using \ac{APF} functions i.e., functions that mimic specific potential field. Considering $G(\delta)$ as a \ac{APF}, it may comprise two components, the attraction $G_a(\cdot)$, and the repulsion $G_r(\cdot)$, and can only converge to a single equilibrium point, occurring at the minimum potential where $G_a(\cdot)=G_r(\cdot)$. Such approach defines controllers laws dependent on a resultant force $f(\cdot)$ due to the potential $G(\cdot)$: $f(\delta) = - \nabla G(\delta)$. This method was already used for trajectory path planning, where a \ac{APF} is inputted in the kinetic models \cite{nag2013behaviour,zhao2020udwadia,d2023distributed,jia2014leader,coquet2021control,ghaderi2023quadrotor}, and also as control laws for dynamic models \cite{kucherov2017group,byun2022potential}. 

\subsection{APF in control}

Eleven \ac{APF} control laws were already proposed, as presented in \autoref{table:APF}.
All the \ac{APF} were applied in the kinematics field, in particular in robotics, in which most of them (7 in 11) are focused on the robot swarm control, and some (4 in 11) are applications related to aircraft trajectory planning.
Considering the attraction behavior, the majority of studies (X in 11) use the quadratic attractor $G_a(\delta)=k\left \| \delta^2 \right \|$, with $k=\frac{1}{2}\lambda_1$ or a similar function only differing in the constant parameter $k$. This formulation was most likely defined due to its simplicity and linearization of the resultant force $f(\delta) = -\nabla G(\delta) = -\lambda_1 \delta$. Regarding the repulsive case, the most common established  function was $G_r(\delta)=\frac{1}{2}\lambda_2\lambda_3  \exp \left({-\left \| \delta^2 \right \| / \lambda_3}\right)$. Even though control dynamics are performed according to attraction and repulsion algorithms, it is important to emphasize that the proposed functions/ control laws are not inspired by potentials with physical or natural significance.

\begin{table}[h!]
\centering
\begin{tabular}{l >{\arraybackslash}p{3cm} >{\arraybackslash}p{4cm} c} 
 \hline
 \textbf{APF Equation} & \textbf{Type} & \textbf{Application} & \textbf{Reference}\\ [0.5ex] 
 \hline
 \rule{0pt}{5ex}   
 $G(\delta) = \frac{1}{2} \lambda_1 \left \| \delta^2 \right \| + \frac{\lambda_2 \lambda_3}{2} e^{-\frac{\left\| \delta \right\|^2}{\lambda_3}} $  & \small Attraction and Repulsion & \small Satellite cluster formation avoiding self collisions  & \cite{nag2013behaviour}  \\
 $G(\delta) = \frac{1}{2} \lambda_1 \left\| \delta \right\|^2 $ & \small Attraction & \small Swarm mechanical system following a target trajectory & \cite{zhao2020udwadia}\\
 $G(\delta) = \frac{1}{2} \lambda_1 \left \| \delta^2 \right \| + \frac{\lambda_2 \lambda_3}{2} e^{-\frac{\left\| \delta \right\|^2}{\lambda_3}} $ & \small Attraction and Repulsion & \small Swarm mechanical system formation avoiding self collisions & \cite{zhao2020udwadia} \\
 $G(\delta) = \frac{\lambda_1}{\left\| \delta \right\|^2} - \lambda_2 \ln{\left( \lambda_3 - \left\| \delta \right\|^2 \right)} + \lambda_4 $ & \small Attraction and Repulsion & \small  Robotic fish leader-follower formation flocking problem & \cite{jia2014leader} \\
 $G(\delta) = \frac{\lambda_1}{2} \left\| \delta_1 \right\|^2 + \frac{\lambda_2}{2} \left\| \delta_2 \right\|^2 + \frac{\lambda_3 \lambda_4}{2} e^{-\frac{\left\| \delta_3 \right\|^2}{\lambda_3}}$ & \small Attraction and Repulsion & \small Path following control of a wireless sensor network avoiding self collision & \cite{byun2022potential} \\
 $G (\delta) = \frac{1}{2} \lambda_1 \left\| \delta \right\|^2 $ & \small Attraction & \small Quad-rotor path control towards a target & \cite{ghaderi2023quadrotor}\\  
 $ G(\delta) = \frac{1}{2} \lambda_1 \left( \frac{1}{\left\| \delta \right\|} -\frac{1}{\lambda_2}\right)^2 \textup{, if} \left\| \delta \right\| \leq \lambda_2  $ & \small Repulsion & \small Quad-rotor obstacles avoidance & \cite{ghaderi2023quadrotor}\\
 $G(\delta) = \lambda_1 \left\| \delta^2 \right \| - \lambda_2 \ln \left( \left \| \delta \right \| \right) $ & \small Attraction and Repulsion & \small Swarm formation for mobile odor source localization problem avoiding self collisions & \cite{coquet2021control}\\
 $ G(\delta) = \lambda_1 \left \| \delta^2 \right \| - \lambda_2 \ln \left( \frac{1}{\lambda_3} \left \| \delta^2 \right \| -1 \right)$ & \small Attraction & \small Agent swarm formation shape control& \cite{d2023distributed}\\
 $ G(\delta) = -\frac{\lambda_1}{2}  \ln \left( \frac{\left \| \delta^2 \right \|}{\lambda_2^2} \right) $& \small Repulsion & \small Agent swarm self collision avoidance& \cite{d2023distributed}\\
 $ G(\delta) = \frac{\lambda_1}{2} \left \| \delta^2 \right \| + \frac{\lambda_2}{2} \left( \frac{1}{\delta}  - \frac{1}{\lambda_3} \right)^2 $& \small Attraction and Repulsion & \small  Movement control of multi-UAV system with leader following and fixed obstacle avoidance & \cite{kucherov2017group}\\[1ex] 
 \hline
\end{tabular}
\caption{APFs focused on control system applications. $G$ is the potential function and $\delta$ the distance between agents and targets / obstacles. $\left \| \cdot \right \|$ is the euclidean norm.}        
\label{table:APF}
\end{table}

\begin{figure}[ht!]
\captionsetup{singlelinecheck = false, justification=justified}
\centering
\sidesubfloat[]{\label{subfig:coquetControl}
\includegraphics[width=0.95\linewidth]{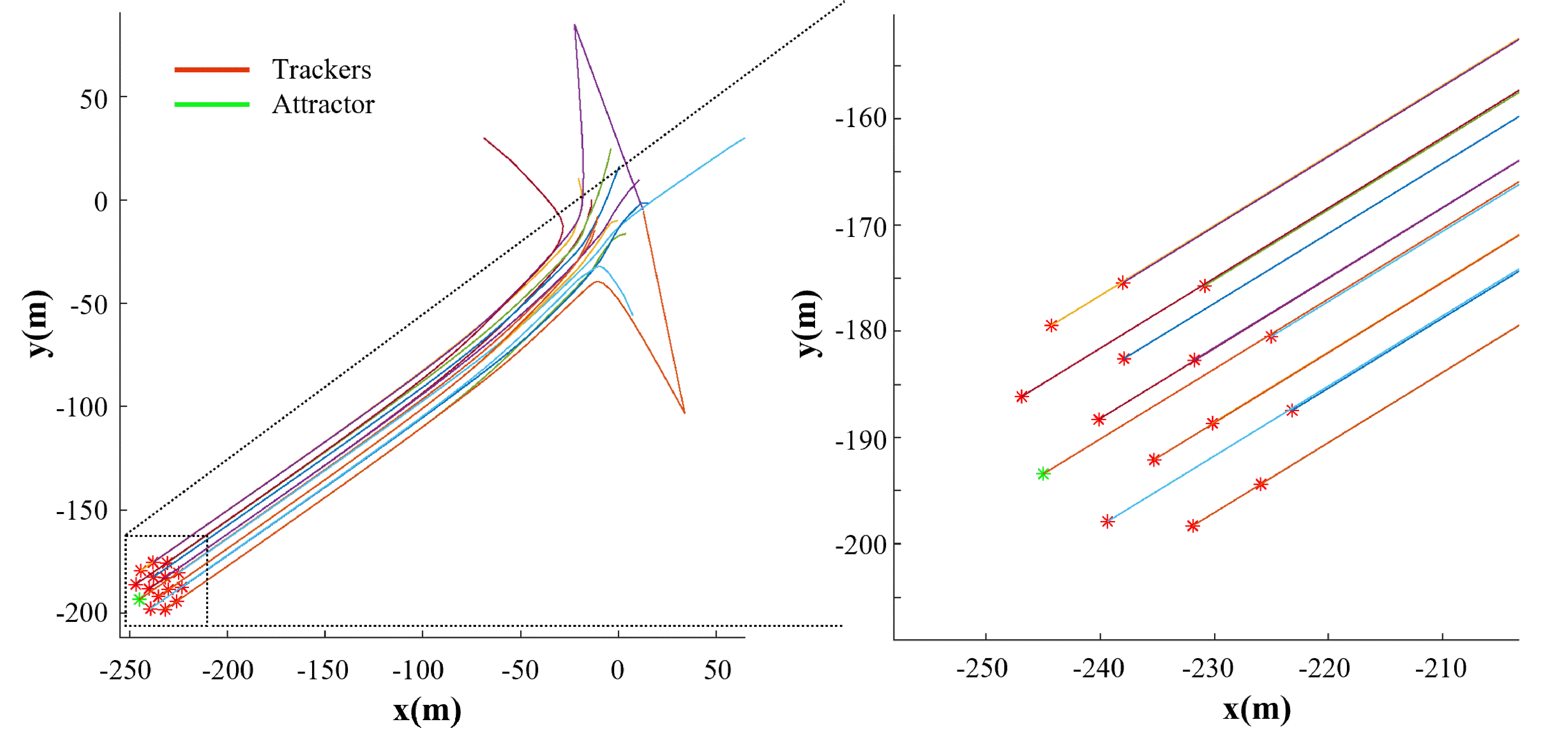}}
\\[20pt]
\sidesubfloat[]{\label{subfig:dalfonsoDistributed}
\includegraphics[width=0.80\linewidth]{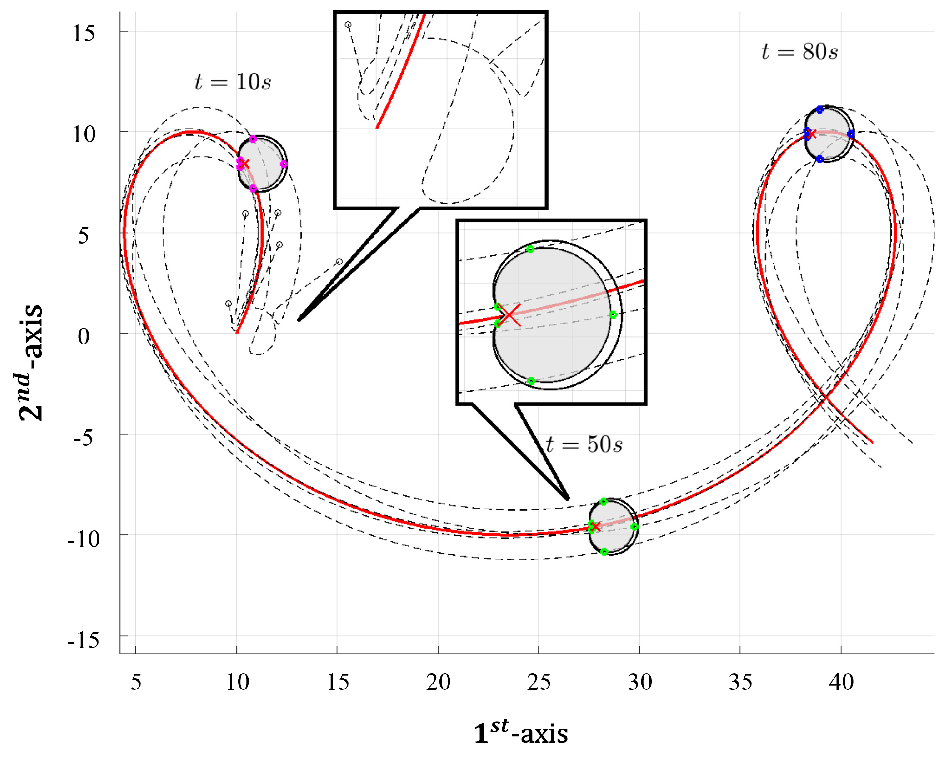}}
\caption{(a) Example of swarm trajectory following the leader (green mark) by using APF control.  Adapted with permission from Ref. \cite{coquet2021control}. (b) Simulation of swarm following a trajectory $\phi(t)$ (red line) by using APF control. The black circles are the initial positions of the agents and the dotted lines are the paths traveled by the agents. Three temporal snapshots of the agents’ states are depicted - magenta: $t = 10~s$; green: $t = 50~s$; blue: $t = 80~s$. Adapted with permission from Ref. \cite{d2023distributed}.}
\label{fig:APF}
\end{figure}

The \ac{APF} was used to control a satellite cluster using the repulsive potential function \sloppy $G_r(\delta)=\frac{1}{2}\lambda_2\lambda_3  \exp \left({-\left \| \delta^2 \right \| / \lambda_3}\right)$, where $\left \| \delta^2 \right \|$ is the distance between two agents \cite{nag2013behaviour}. Self-collisions are avoided by summing the resulting potentials between all pairs of agents. A similar approach was used to control a mobile robot swarm to track moving target \cite{zhao2020udwadia}, as well as in a multiple fish robot system with a leader \cite{jia2014leader}. The capability of APF control to guide a swarm along a designated trajectory and ensure formations with specific shapes is illustrated in \autoref{subfig:coquetControl} and \autoref{subfig:dalfonsoDistributed}. The attractive function $G_a(\delta)=\frac{1}{2}\lambda_1\left \| \delta^2 \right \|$ causes the agents to attract each other, generating the formation region. The swarm region is unique as the controlled systems are moved towards an equilibrium state of minimum potential. Considering the ideal kinetic model, the application of \ac{APF} is given by $v(t+1) = v(t) - \nabla G(\delta)$, where $- \nabla G(\delta)$ defines the acceleration established by the dynamics related to the kinematic model \cite{zhao2020udwadia,jia2014leader,coquet2021control}. \pagebreak

\captionsetup[figure]{justification = centering}
\captionsetup[subfigure]{justification = centering}

In such cases, the \ac{APF} method is used to generate desired trajectory, but a dynamic controller is needed to effectively track the defined trajectories ensuring robustness, i.e., a controller with compensating ability in response to disturbances deviating the controlled system from desired trajectories. Differently, the \ac{APF} was designed taken $G_a(\delta)=- \lambda_2 \ln{\left( \lambda_3 - \left\| \delta \right\|^2 \right)}$ and $G_r(\delta)=\frac{\lambda_1}{\left\| \delta \right\|^2}$, where $\left\| \delta \right\|$ is also the distance between agents \cite{jia2014leader}. In such case, a \ac{FLC} was designed for path tracking defined by the potential function method. When the \ac{APF} is inputted to the dynamic model, the gradient of the potential function is introduced on the control input of the dynamic model, as follows: $x(k+1) = f(x(k)) + u(k); u(k+1) = u(k) - \nabla G(\delta)$ \cite{kucherov2017group,byun2022potential}, where $x(k)$ and $u(k)$ are the current state and the control input, respectively. 
The \ac{APF} was also incorporated in the dynamic model, aiming to control mobile collectors \cite{byun2022potential}. The proposed function is composed by three components, two of them implementing attractiveness dynamics, and another one implementing repulsiveness dynamics. 
The first attractive component is used to eliminate the state estimation error, by defining $\frac{\lambda_1}{2} \left\| \delta_1 \right\|^2$, where $\delta_1$ is given by $x(k) - \hat{x}(k)$, and $\hat{x}(k)$ is the estimated state.
The second component was defined as $\frac{\lambda_2}{2} \left\| \delta_2 \right\|^2$, with $\delta_2 = x(k) - x_d(k)$, where $x_d(k)$ is the reference or the desired system state, such that it is able to move the system under control towards the desired states. 
In the repulsive component, $\delta_3$ is the distance between the agents, and it is used to avoid self-collisions. The strategy based on incorporating the \ac{APF} in the control law is highly promising to tackle feedback systems, involving the fundamental concept of attraction of system state towards the required state according to a natural rationally found in the real physical universe.


\section{Discussion}

\begin{figure}[ht!]
\captionsetup{singlelinecheck = false, justification=justified}
\centering
\sidesubfloat[]{\label{subfig:pieGSA}\includegraphics[width=0.43\linewidth]{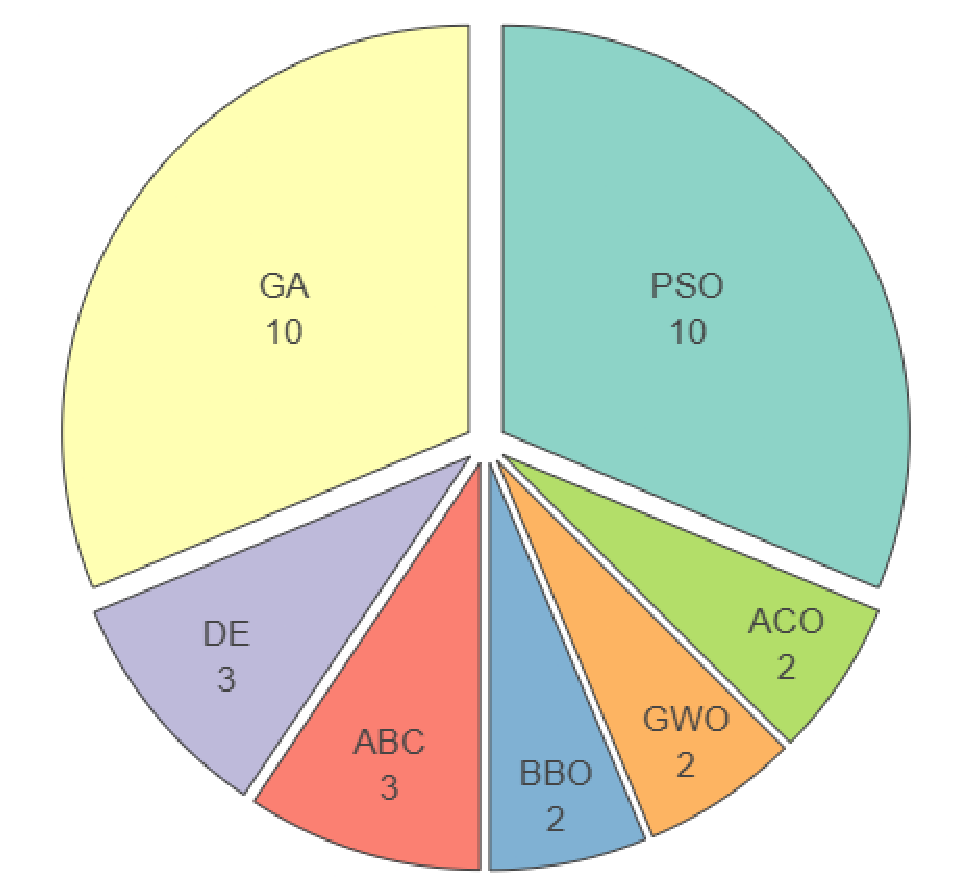}}
\sidesubfloat[]{\label{subfig:gaoDistributed}\includegraphics[width=0.43\linewidth]{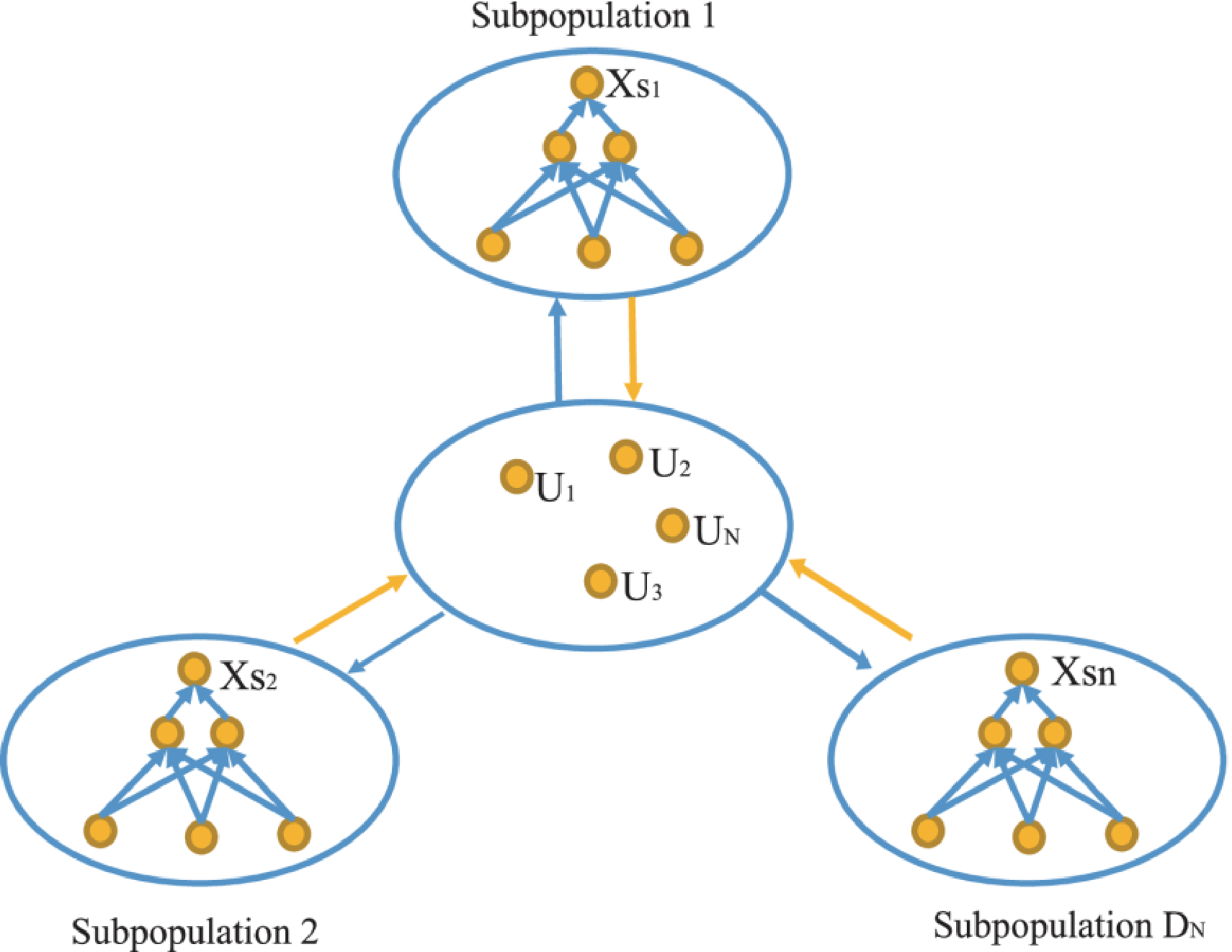}} \\ [5ex]
\sidesubfloat[]{\label{subfig:wantMultilayer}\includegraphics[width=0.58\linewidth]{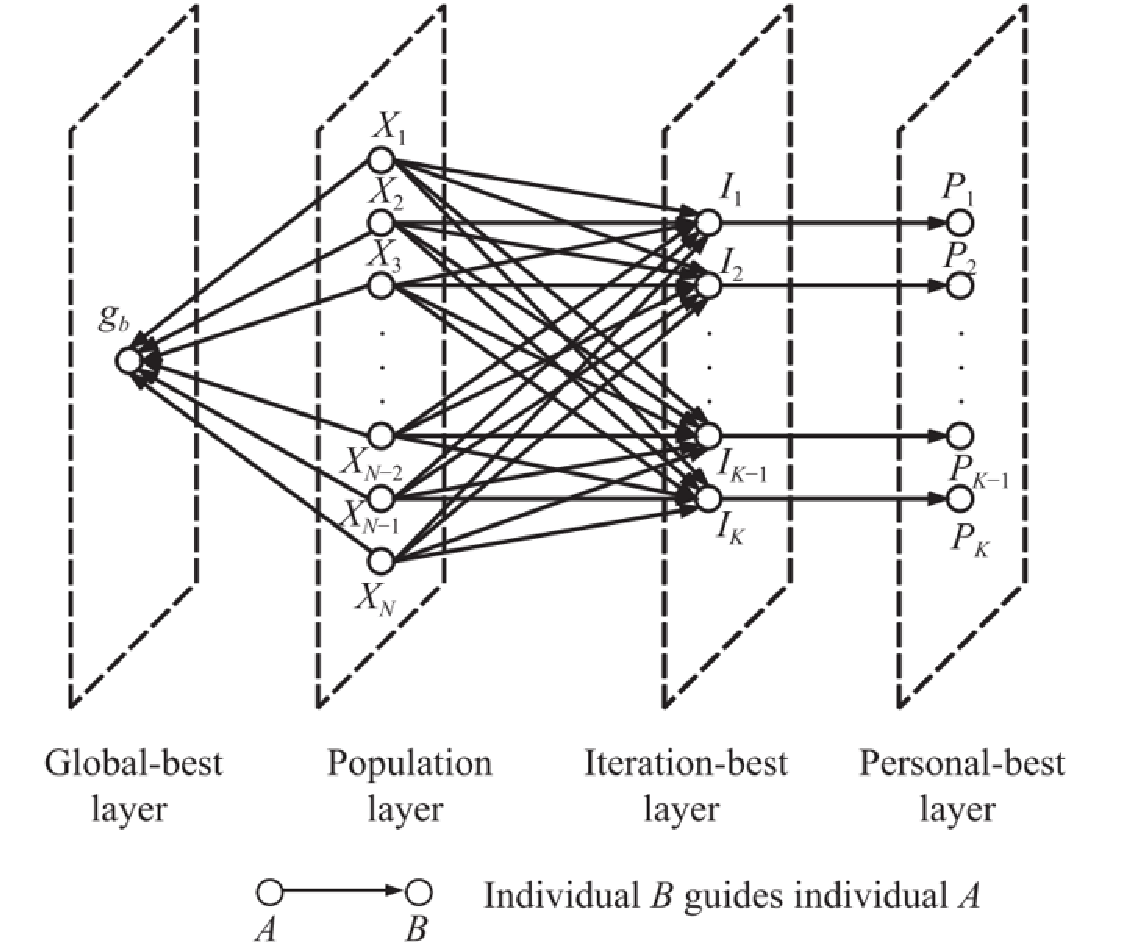}}
\caption{(a) The number of studies in the control field wherein the \ac{GSA} exhibits a superior optimization capacity compared to other algorithms documented in the literature. (b) Population distributed structure of Distributed Multi-Layer Gravitational Search Algorithm. Reproduced with permission from Ref. \cite{li2021novel}. (c) Illustrative population structure of Hierarchical Multi-Layered Gravitational Search Algorithm. Reproduced with permission from Ref. \cite{wang2020multi}.}
\label{fig:discussion}
\end{figure}

Most studies perform comparisons with other well-established algorithms to demonstrate their performance. Nonetheless, given the wide range of existing metaheuristics and optimization methods, conducting a comprehensive comparison with all of them is hardly feasible, making it difficult to conclusively determine which is the best algorithm. Based on the analyses here performed, the GSA exhibited superior performances in 20 (out of 20) control problems  when compared to \ac{PSO} and \ac{GA}. Additionally, it outperformed \ac{DE} and \ac{ABC} in 6 (out of 6) control problems. These findings, depicted in \autoref{subfig:pieGSA},  suggest the superiority of \ac{GSA} over \ac{PSO} and \ac{GA}. While other comparisons are presented, they lack significant expression for a meaningful analysis. Indeed, the \ac{BH} demonstrated a performance superiority over \ac{PSO} in 2 cases (out of 2), and \ac{MVO} outperformed \ac{DE}, \ac{PSO}, \ac{GA}, and \ac{GWO} in 4 cases (out of 4). This discrepancy highlights the disparity in the number of studies applied to control between \ac{GSA} and other methods inspired by gravitational attraction (\ac{BH} and \ac{MVO}). Even so, it is important to recall the No Free Lunch theorem \cite{wolpert1997no}, which state that  there is not a single algorithm that outperforms all others across all metrics and problems:  performance is inherently problem-dependent. 

The GSA provides several advantages, namely its minimal need for hyper-parameter adjustment, as it only requires three parameters: $G_0$, $\beta$, and $\varepsilon$. Moreover, it typically reaches rapid convergences  even though its non-complex implementation (see \autoref{fig:gsaimproved}d), and features  intuitive and easily interpretable movements of agents. GSA also facilitates the incorporation of adaptation mechanisms, as well as the resolution of constrained problems. However, it can exhibit low agent dispersion during the final stages, resulting in reduced exploration capabilities. Additionally, GSA lacks memory, and may become trapped in local minima during the latter iterations. Other limitations include the blind evolution of the variable $G$, and the need to predefine the maximum number of iterations, as the evolution of $G$ is dependent of the duration of the established optimization process. 

Concerning the BH algorithm, it has been formulated using a straightforward mechanism, easily understandable  and implemented, only requiring the computation of two equations. Notably, the BH is devoid of hyper-parameters, ensuring performances independent of user settings. This algorithm effectively tackles problems with constraints, and maintains a strong exploration capacity across all iterations. When a solution is eliminated, it is promptly replaced by a new one randomly generated within the search space, ensuring population dispersion. Despite its simplicity, its few parameters hardly allow the inclusion of adaptation mechanisms in its original formulation, as well as the ability to effectively balance exploration and exploitation. It relies solely on randomization, potentially leading to extended convergence times. Furthermore, its formulation lacks is not supported by real physical laws, which significantly deviates from an authentic alignment with natural inspired rationality.

The MVO is recognized for its rapid convergence and robust exploration capabilities, facilitated by solution mutation and crossover. It only requires one hyper-parameter to determine its exploration ability. However, MVO poses relevant challenges in interpretation, due to its complex nature and hypothetical astronomical objects. Moreover, it shares the limitation of requiring the predefinition of the maximum number of iterations. Additionally, auxiliary mechanisms  need to be incorporated into its algorithm, such as the roulette wheel selection and sorting procedures, the latter being computationally heavier due to the sorting of solutions in each iteration.

The novelty of \ac{GSO} is indeed questionable, given that its architecture bears resemblance to multi-layered versions of other algorithms. Additionally, its classification as an optimization algorithm is debatable, as the proposed concept lacks an inherent optimization mechanism and relies on other mechanisms. Therefore, \ac{GSO} can be perceived more as an advancement over existing algorithms, rather than a novel optimization approach.
Furthermore, while \ac{GSO} claims inspiration from gravitational attraction between stars and galaxies, its attraction mechanisms are contingent upon the underlying algorithm, which may or may not incorporate gravitational principles. This reliance on a external algorithms limit the argumentation presented by the authors presenting the \ac{GSO} as an optimization method.
Although the potential of utilizing \ac{GSO} to enhance the performance of algorithms like \ac{GSA} or \ac{BH} remains unexplored, such an endeavor holds strong motivation for future research. These algorithms demonstrate complementary foundational inspirations, suggesting that exploring their integration could yield promising results.

To address the mentioned limitations, researchers have been developing enhancements and modifications, introducing new mechanisms and hybridization, as highlighted in \autoref{table:variationsGSA} and \autoref{table:BH}. However, concerning variants of algorithms applied to control, the most relevant developments have been focused on the use of GSA. In the case of BH algorithms, only two variants have been identified, while only the original MVO algorithm was already applied in control applications. Therefore, it is mandatory to highlight other variants aiming to implement more sophisticated mechanisms, even if they remain untested in the field of system control, as they hold potential to stimulate future excellence research in a multidisciplinary basis.  Relevant studies are detailed in \autoref{table:noControlAlg}.
Included are other \ac{GSA} relevant variants: (1) a Curved Space \ac{GSA}, where a dimension reduction technique is applied to the problem search space \cite{Giladi2020}; (2) a Memory-based \ac{GSA}, where the concept of personal best is introduced \cite{Darzi2016}; (3) adaptive versions, that introduce novel mechanisms to guide agents out of local optimum trapping \cite{guo2023adaptive}, and a personal gravitational constant \cite{lei2020aggregative}; (4) hierarchical and distributed \ac{GSA}, where the overall population is divided in subsets of populations (\autoref{subfig:gaoDistributed}) with multi layers of hierarchy (\autoref{subfig:wantMultilayer}) \cite{wang2019hierarchical,wang2020multi,li2021novel,wang2021gravitational}; and (5) Multiple chaos \ac{GSA} a mechanism that uses several chaotic maps \cite{song2017multiple,gao2014gravitational}, supported by the performance superiority of multiple embedded chaos \cite{gao2019chaotic}. Regarding \ac{BH}, other variants include: (1) an improved \ac{BH}, which includes a crossover mechanism, inspired in Genetic algorithms, to generate new agents, avoiding a random generation \cite{deeb2022improved}; (2) a Chaotic inertia weight \ac{BH} which improves the local search by using chaotic maps \cite{pashaei2021medical}; and (3) a Multi-population \ac{BH} which uses multiple populations of agents, instead of a single one \cite{salih2023multi}. To enhance \ac{MVO}, it was already proposed: (1) Chaotic \ac{MVO} \cite{ewees2019chaotic}, by introducing the chaotic behavior in the standard \ac{MVO}, as already proposed in other algorithms \cite{yuan2014hybrid,gandomi2013firefly,alatas2009chaos} ; and (2) a hybrid Sine-Cosine-\ac{MVO} algorithm \cite{otair2022optimized,jui2020modified}.

Given the wide range of metaheuristics and their variants \cite{thymianis2022integration,tzanetos2021nature}, the question arises: is modification and hybridization the path to higher-performance algorithms? Thyamianis \textit{et al.} \cite{thymianis2022integration} reported evidence suggesting that hybridization and the inclusion of additional mechanisms can have positive effects on nature-inspired algorithms. However, they highlight that the additional algorithm complexification can also result in not relevant improvements in exploration and exploitation. Piotrowski \textit{et al.} \cite{piotrowski2018some} also question the complexity inherent to innovative hybrid algorithms: as complexity of the numerous modifications to basic algorithms increases, the risk of  discouraging their use also increases. There is supporting evidence suggesting that simplifying algorithms can enhance transparency and performance \cite{piotrowski2018some}.Indeed, modifications to nature-inspired algorithms often introduce artificial (not-natural inspired) elements into basic algorithms, in a clear opposition to their primary foundations, which stemmed from their assumed simplicity of interpretation and use, grounded in natural origins. Another relevant contradiction identified in the analyzed algorithms is their formulation without explicit use of real astrophysics laws, even though they are labeled as nature-inspired algorithms. Notice that \ac{BH} and \ac{MVO} algorithms only incorporate nature-inspired concepts without a rigorous mathematical foundation of the related phenomena. Concerning the GSA, although it seems to have the strongest connection to the real astrophysical nature, Gauci \textit{et al.} \cite{gauci2012gsa} concluded that it cannot be truly inspired by Newton's laws of gravity, because the square of the distance is disregarded. Indeed, there is strong evidence suggesting that the force model formulated so far for the \ac{GSA} algorithms does not rely on the distance between agents at all \cite{gauci2012gsa}. Thus, the movement mechanism in the \ac{GSA} is primarily proportional to the fitness of the solutions, as the division by the distance between agents mainly serves to normalize it. Based on this analysis, there is evidence indicating that the movement of agents in the \ac{GSA} bears similarities to that of the \ac{PSO}.

Finally, we must highlight that, the use of APF emerged as the most closely related approach to leveraging gravitational attraction or the dynamics of black holes in the development and design of controllers. The results obtained through the combination of attraction and repulsion functions are intuitive and straightforward to interpret, rendering this approach particularly appealing. Nevertheless, the use of dynamics inspired in Newtonian or Einsteinian gravitation, including as it occurs in black holes, in control systems has not yet been deeply researched: they have only been an inspiration-trigger  only aiming to design new methods for optimization problems of well-known control methods (\textit{e.g.} \ac{PID} and \ac{FLC} controllers). To our knowledge, no control methods were developed whose formalization is directly and truly inspired by gravitational attraction laws, and related mathematical Newtonian/ Einsteinian-based formulation. 

\begin{landscape}
\thispagestyle{empty}
\begin{longtable}{>{\footnotesize\arraybackslash}p{1cm} >{\footnotesize\arraybackslash}p{2cm} >{\footnotesize\arraybackslash}p{16.5cm} c} 
\caption{Other recent relevant variants of \ac{GSA}, \ac{BH} and \ac{MVO}, not applied in the control field.}        
\label{table:noControlAlg}
\\
\hline
\textbf{Base method} & \textbf{Designation} & \textbf{Concept novelty description} &\footnotesize\textbf{Reference}\\[0.5ex]
\hline
\ac{GSA} & Curved Space \ac{GSA} & Feasible solutions of \ac{GSA} may be contained in a manifold of lower dimensions, in such a way that, according to the Euclidean distance, the agents may seem close, although they are far apart. The proposed modification calculates the distance between agents along the manifold, instead of directly calculate the Euclidean distance, by utilizing diffusion maps as dimensionality reduction & \cite{Giladi2020} \\
\cline{2-4}
 & Memory-based \ac{GSA} & This algorithm ensures that the best position of any agent ($pbest$) is stored as the agent’s personal best position, and thus, the new positions of the agents, are always calculated based on the previous best values, such that the path towards the best solution is not lost. In this formulation $R_{ij}(t)$ is computed by: $R_{ij}(t) = \left \| X_i(t),pbest_j(t) \right \|$ & \cite{Darzi2016} \\
\cline{2-4}
 &  Adaptive position-guided \ac{GSA} (disGSA) & A novel mechanism is proposed to guide agents out of local optimum trapping in the direction of global best solution: $v_i^d(t+1) = rand \times v_i^d(t) + c_1 a_i^d(t) + c_2 a_{i2}^d(t)$, where $a_{i2}^d(t) = \frac{F_{i2}^d(t)}{M_i(t)}$ and $F_{i2}^d(t) = \sum_{j\in Dbest,j\neq i }^{} rand_j F_{ij}^d(t)$. $Dbest$ is a sorted set of $D_i$ defined as $D_i(t) = \frac{R_{i,best}(t)}{R_{i,worst}(t)+\varepsilon}$. &  \cite{guo2023adaptive} \\
\cline{2-4}
 & Self-adaptive and aggregative learning \ac{GSA} & 
This method proposes an adaptive mechanism wherein each agent possesses its own gravitational constant, defined as: $G_i(t) = \left\{\begin{matrix} G_i(t)r_i(t) & \textup{if}~counter > \theta~\textup{and}~rand<p \\ G_i(t) & \textup{otherwise}\end{matrix}\right.$, where $r_i(t) = \left | \textup{log}\left ( \frac{\left | a_i(t) \right |}{G_i(t)} \right ) \right |$. An aggregative mechanism was also included, with (\ref{eq_totalForce}) being replaced by: $F_i^d(t) = \sum_{j=1}^{k}Y_i(j)$, with $Y_i(j) = \frac{\sum_{k=1}^{j}G_k(t)}{j}\left \{ \frac{M_i(t)}{R_{i,k}+\varepsilon}\sum_{k=1,k\neq j}^{j}r_k M_k(t)\left [ x_k^d(t)-x_i^d(t) \right ] \right \}$ & \cite{lei2020aggregative} \\
\cline{2-4}
 & Multi hierarchical layer \ac{GSA} & 
A hierarchical population structure categorizes individuals into layers based on specific criteria, guiding their evolution in a systematically basis. Layers, organized from top to bottom like a tree, influence individuals progressively. This hierarchical arrangement fosters interactive relationships among layers, shaping the evolution of the population. & \cite{wang2019hierarchical,wang2020multi} \\
\cline{2-4}
 & Distributed and hierarchical \ac{GSA}&  In addition to the hierarchical layers with varying levels, the population is also distributed into several sub-populations, each one segmented into hierarchical layers. & \cite{li2021novel,wang2021gravitational}\\
\cline{2-4}
 & Multiple chaos \ac{GSA} & Multiple chaotic maps into the \ac{GSA} are incorporated as follows: (1) In each iteration, a new chaotic map is selected randomly; (2) The agents undergo mutation using all the different chaotic maps under consideration. Among the solutions generated, the one with the best fitness value is preserved; (3) The probability of selecting a specific chaotic map in each iteration is dynamically adapted based on its success rate. & \cite{song2017multiple,gao2014gravitational}\\
\hline
\ac{BH} & Improved \ac{BH} & With a certain probability, when a new agent is created, the process involves crossing over two existing feasible solutions instead of randomly generating a new one. & \cite{deeb2022improved} \\
\cline{2-4}
 & Chaotic inertia weight \ac{BH} & This algorithm uses chaos theory to enrich the search behavior. A hiper-parameter named inertia weight ($w$) is introduced to control the balance between exploration and exploitation, with $w$ given by: $\text{w}(t) = (\text{w}_{max} - \text{w}_{min}) \left (\frac{t_{max}-t}{t_{max}} \right) +\text{w}_{min}C(t)$, where $C(t)$ is a chaotic map. Therefore, (\ref{eq_pos_bh}) is replaced by $x_i(t+1)=\text{w}(t) x_i(t)+rand \; \left( x_{BH} - x_i(t) \right)$  & \cite{pashaei2021medical} \\
\cline{2-4}
 & Multi-Population \ac{BH} & This algorithm has the same formulation of the base version but uses multiple populations instead of a single one. At the end, the solution is the best agent of all populations & \cite{salih2023multi} \\
\hline
\ac{MVO} & Chaotic \ac{MVO} & This algorithm proposes to replace $r_4$ in (\ref{eq:r4mov}) by a chaotic map, to improve the local search ability of standard \ac{MVO} & \cite{ewees2019chaotic} \\
\cline{2-4}
& Sine Cosine \ac{MVO} & By using a sine cosine mechanism, (\ref{eq:r4mov}) is replaced by: \newline
 $x_i^j = \left\{\begin{matrix}
\left\{\begin{matrix}
\textup{AP}+ \textup{TDR}\cdot \left( \left| 2r_6 X_j -x_i^j \right|\text{sin}(2\pi r_5) \right) & r_3 <0.5 \\ 
\textup{AP}- \textup{TDR}\cdot \left( \left| 2r_6 X_j -x_i^j \right|\text{cos}(2\pi r_5) \right) & r_3 \geq 0.5
\end{matrix}\right. & r_2< \textup{WEP}\\ 
x_i^j & r_2\geq \textup{WEP}
\end{matrix}\right.$, where $\textup{AP} = \left( X_j + x_i^j \right)/2$ \newline
& \cite{otair2022optimized,jui2020modified} \\ 
\hline
\end{longtable}
\end{landscape}

\section{Conclusions}

Significant scientific breakthroughs have been carried out in the field of Control Engineering using universe-inspired algorithms. Two main categories have been the focus of such advances: optimization algorithms applied in control problems (\ac{GSA}, \ac{BH}, \ac{MVO}, and \ac{GSO}), where main improvements were achieved in the scope of control parameters optimization; and the identification of a proper control method, inspired by the attraction between bodies, known as Artificial Potential Fields,  which was mainly used to guide agents to an equilibrium state defined by choosing appropriate attraction and repulsion functions.

\ac{GSA} algorithms has been designed according to the  gravity law, and the movement of the agents is due to gravitational forces, which allow the information transfer between agents, as masses within the gravitational system are affected by one another.  Most results obtained by GSA were able to provide superior results in comparison with \ac{GA} and \ac{PSO}.  Concerning the \ac{BH} algorithms, two significant advantages were identified: (i) its structure is not complex, and its  implementation is easily performed; (ii) it does not raise parameter tuning issues. For these reasons, it is considered a feasible option  when fast and accurate results are required within a short period of time. \ac{MVO} and \ac{GSO} are the most recent proposed algorithms, and, for this reason, the number of studies analyzing their performance and characteristics remains limited. Applications and modifications carried out to GSA, BH, MVO, and GSO algorithms have revealed that such concepts can be used in a wide range of control problems, and can still evolve towards improved performances. APF represents the closest approach to natural gravitational phenomena by introducing artificial attraction behaviors.  Therefore,  future high-sophisticated control systems inspired by black-hole attraction dynamics can be engineered if they are further analyzed and considered.

Some difficulties arise when effective conclusions must be stated. On the one hand, only few studies present the data in a clear way or report meaningful comparisons; on the other hand, relevant data is lacking in most studies (such as convergence times). Besides, the comparison between the different approaches here analyzed is hard to achieve, due to the influence of: (i) the diversity of methodologies ; (ii) the diversity of the objective functions; (iii) the parametrisation of the algorithms, as different parameters can conduct to different results; and (iv) the diversity of the applications and related scopes. Despite all these problems,  the achieved results suggest that optimization and proper control methods inspired by gravitation and black holes attraction perform better than other approaches, namely the \ac{PSO}, while ensuring easier implementation and interpretation. Nevertheless, these results highlight the capability of black holes, gravitational attraction, and universe dynamics in general, to overcome many control engineering problems, even though they are still limited to the field of optimization and metaheuristics. Likewise, future studies may explore realist universe-related dynamics in order to design effective control methods. To date, no control methods have been truly designed  from scientific formulations related to real astrophysical phenomena. This fact, together with the results achieved in this study, provides new research directions where  highly innovative concepts can be developed, namely controller ruled by astrophysics-like laws to establish effective bridges between black hole physics and automatic control. Attraction may behave as a feedback mechanism of the distance between the considered masses. Consequently,  future controllers can be built upon natural feedback interactions inspired by the gravitational attraction towards the singularity of black holes. Supported by analogies with physics, where nothing can escape from a black hole once the Schwarzschild radius is crossed, new concepts can established  such stable equilibrium points.. These are highly promising future prospects that overcome the methods here discussed and analyzed,  as they do not employ  control approaches based on artificial phenomena, avoiding then to neglect the rationality and effectiveness inherent of natural systems.

\FloatBarrier


\section*{Data availability statement}

No data was used for the research described in the article.


\section*{Declaration of competing interest}

The authors declare that they have no known competing financial 
interests or personal relationships that could have appeared to 
influence the work reported in this paper.


\section*{Acknowledgments}

R.M.C.B. and M.P.S.S. were supported by the Portuguese Foundation for Science and Technology (FCT) (PhD grant reference: 2023.01947.BD; project references: UIDB/00481/2020 and UIDP/00481/2020; DOI 10.54499/UIDB/00481/2020, https://doi.org/10.54499/UIDB/00481/2020, and DOI 10.54499/UIDP/00481/2020, https://doi.org/10.54499/UIDP/00481/2020). D.F.M.T. and C.A.R.H were supported by FCT through project UIDB/04106/2020 (DOI 10.54499/UIDB/04106/2020, https://doi.org/10.54499/UIDB/04106/2020).



\end{document}